\documentclass[final,onefignum,onetabnum]{siamart220329}

\usepackage{lipsum}
\usepackage{amsfonts}
\usepackage{graphicx}
\usepackage{epstopdf}
\usepackage{algorithmic}
\usepackage{amsmath}
\usepackage{dutchcal} 
\usepackage{amsopn}

\ifpdf
\DeclareGraphicsExtensions{.eps,.pdf,.png,.jpg}
\else
\DeclareGraphicsExtensions{.eps}
\fi

\usepackage{xcolor}
\definecolor{darkgreen}{rgb}{0.1,0.6,0.1}

\newcommand{\vv}{\mathbf{v}}  
\newcommand{\ee}{E} 
\renewcommand{\S}{S} 
\newcommand{\temp}{\theta} 
\newcommand{\mde}[2]{\frac {\text{d} #1}{\text{d}#2}}
\newcommand{\halb}{\frac{1}{2}}
\newcommand{\x}{\mathbf{x}}  
\newcommand{\q}{\mathbf{u}} 
\newcommand{\w}{\mathbf{w}} 
\newcommand{\n}{\mathbf{n}} 
\newcommand{\dt}{\Delta t}
\newcommand{\lnpc}{l_{pc} \n_{pc}} 
\newcommand{\lnpb}{l_{pb} \n_{pb}} 
\newcommand{\Mpc}{\mathbf{M}_{pc}} 
\newcommand{\Mp}{\mathbf{M}_{p}} 

\allowdisplaybreaks

\newsiamremark{remark}{Remark}
\newsiamremark{hypothesis}{Hypothesis}
\crefname{hypothesis}{Hypothesis}{Hypotheses}
\newsiamthm{claim}{Claim}

\headers{A new HTC scheme for Lagrangian gas dynamics}{W. Boscheri and M. Dumbser and P.H. Maire}    %

\title{A new thermodynamically compatible finite volume scheme for Lagrangian gas dynamics}

\author{Walter Boscheri \thanks{Department of Mathematics and Computer Science, University of Ferrara, via Machiavelli 30, 44121 Ferrara, Italy (\email{walter.boscheri@unife.it})} 
	\and 	
	Michael Dumbser\thanks{Department of Civil, Environmental and Mechanical Engineering, University of Trento, Via Mesiano 77, 38123 Trento, Italy (\email{michael.dumbser@unitn.it})}
	\and 	
Pierre-Henri Maire \thanks{CEA CESTA, 33116 Le Barp, France (\email{pierre-henri.maire@cea.fr})} 
}

\begin{document}

\maketitle

\begin{abstract}
The equations of Lagrangian gas dynamics fall into the larger class of overdetermined hyperbolic and thermodynamically compatible (HTC) systems of partial differential equations. They satisfy an entropy inequality (second principle of thermodynamics) and conserve total energy (first principle of thermodynamics). 
The aim of this work is to construct a novel thermodynamically compatible cell-centered Lagrangian finite volume scheme on unstructured meshes. Unlike in existing schemes, we choose to directly discretize the entropy inequality, hence obtaining total energy conservation as a consequence of the new thermodynamically compatible discretization of the other equations. First, the governing equations are written in fluctuation form. Next, the non-compatible centered numerical fluxes are corrected according to the approach recently introduced by Abgrall \textit{et al.}, using a scalar correction factor that is defined at the nodes of the grid. This perfectly fits into the formalism of nodal solvers which is typically adopted in cell-centered Lagrangian finite volume methods. Semi-discrete entropy conservative and entropy stable Lagrangian schemes are devised, and they are adequately blended together via a convex combination based on either \textit{a priori} or \textit{a posteriori} detectors of discontinuous solutions. The nonlinear stability in the energy norm is rigorously demonstrated and the new schemes are provably positivity preserving for density and pressure. Furthermore, they exhibit zero numerical diffusion for isentropic flows while being still nonlinearly stable. The new schemes are tested against classical benchmarks for Lagrangian hydrodynamics, assessing their convergence and robustness and comparing their numerical dissipation with classical Lagrangian finite volume methods.   

\end{abstract}



\section{Introduction}\label{sec.intro}
The gas dynamics equations constitute a nonlinear hyperbolic system of conservation laws that evolve mass, momentum and total energy. As discussed in \cite{Lax1972}, the solutions might exhibit discontinuities and shock waves at finite time, which necessitate the introduction of weak solutions.  The physically meaningful solution is then identified relying on an admissibility criterion known as entropy condition. This permits to reject nonphysical solutions by enforcing consistency with the second law of Thermodynamics, thus ensuring entropy conservation for smooth flows and entropy production, which transfers kinetic energy into internal energy, in the case of discontinuous flow fields. These systems are said to be hyperbolic thermodynamically compatible (HTC). It is worth recalling that the connection with thermodynamic compatibility was first presented in \cite{God1961} and later in \cite{FriedrichsLax}, noticing that the gas dynamics equations belong to the more general framework of symmetrizable hyperbolic systems \cite{FriedrichsSymm}. 

To numerically solve the hydrodynamics system, the updated Lagrangian form is considered in this work. The Lagrangian approach allows strong shocks and expansions to be properly modeled for compressible fluid flows due to the embedded mesh motion, and material interfaces or contact waves can be accurately tracked since zero mass flux occurs across the boundaries of the control volumes. This peculiar representation is characterized by a frame moving with the material velocity. The first Lagrangian numerical method dates back to the 50s \cite{Neumann1950}, where a staggered grid was employed for the kinematic variables, hence requiring an artificial viscosity term to be tuned in order to deal with shock waves. An alternative approach is given by Godunov-type finite volume (FV) methods, in which the integral form of the conservation laws is discretized on a collocated grid. As such, numerical fluxes have to be computed at the cell interfaces by means of Riemann solvers \cite{HLL1983}. 
Mimicking the thermodynamic compatibility at the discrete level is not trivial. Staggered grid Lagrangian methods discretize the internal energy equation, hence ensuring entropy conservation for smooth flows up to the order of the numerical method. However, entropy dissipation must be added for discontinuous solutions via an artificial viscosity that needs to be properly adjusted so that thermodynamic compatibility is fulfilled in multiple space dimensions. On the other hand, Godunov-type schemes are compliant with an entropy inequality within each cell, even in the case of isentropic flows, hence yielding severe inaccuracy for strong rarefaction waves \cite{Kulikovskii}. The design of thermodynamically compatible finite volume schemes has been investigated in the Eulerian framework in \cite{Tadmor1}, where the total energy equation is discretized and the entropy equation (or inequality) is enforced by designing specific numerical fluxes at the interfaces. Following this seminal idea, entropy-compatible schemes have been developed in the Eulerian setting, see for instance \cite{ray2016,FjordholmMishraTadmor,Tadmor2003,Hiltebrand2014,GassnerSWE,GassnerEntropyEuler,ShuEntropyMHD1}. 

Very few contributions can be found in the literature concerning the construction of Lagrangian methods that are compliant with the second law of thermodynamics \cite{BRAEUNIG2016127,Maire2020}. Therefore, the aim of this work is to design a new Lagrangian finite volume scheme that is provably thermodynamically compatible. However, unlike most of the existing numerical methods for hyperbolic systems, we do \textit{not} discretize the total energy equation, but we directly discretize the \textit{entropy inequality} and instead obtain total energy conservation as a \textit{consequence} of the novel compatible discretization, thus admitting a discrete extra conservation law to the governing equations. This approach has first been proposed in \cite{SWETurbulence,HTCGPR,HTCMHD} for turbulent shallow water flows \cite{Ivanova2019}, the GPR model of continuum mechanics \cite{PeshRom2014} and the MHD equations, which both also fall into the wider class of hyperbolic and thermodynamically compatible (HTC) systems. A general framework for the construction of schemes that satisfy additional extra conservation laws has been forwarded in \cite{Abgrall2018}, and this technique has been recently also used in \cite{HTCAbgrall} in the context of finite volume and discontinuous Galerkin methods. Following the ideas of Abgrall, in this work we will compute a node-based scalar corrector factor that makes the scheme compatible with thermodynamics in the Lagrangian setting, demonstrating that total energy is retrieved as an additional conservation law. To design entropy stable schemes, we propose to blend the new thermodynamically compatible methods with the EUCCLHYD scheme \cite{EUCCLHYD}, rewritten here in fluctuation form and for the entropy inequality. The production term is then exactly quantified, and the positivity preservation of density and pressure can be proven. We stress that in our framework the entropy production is directly controlled by the scheme, since here the entropy inequality is directly solved for the first time in the Lagrangian context, which is a radically different concept compared to most of the existing entropy compatible numerical methods. We will focus on the design of semi-discrete schemes, thus keeping the time continuous while correcting the spatial fluxes to obtain thermodynamic compatibility.

The rest of this paper is organized as follows. In Section \ref{sec.model} we present the governing equations, starting from the energy equation and passing to the entropy inequality. Section \ref{sec.method} is devoted to the construction of the new hyperbolic thermodynamically compatible (HTC) Lagrangian scheme, which introduces both semi-discrete entropy conservative and entropy stable schemes, including the blending procedure. Nonlinear stability in the energy norm as well as positivity of thermodynamics quantities are rigorously demonstrated. Next, in Section \ref{sec.results} we present some numerical results for a suite of benchmarks for Lagrangian hydrodynamics. Finally, conclusions and an outlook to future investigations are drawn in Section \ref{sec.concl}.     

\section{The Euler equations of Lagrangian hydrodynamics} \label{sec.model}
Let $\x \in \mathbb{R}^d$ be the spatial position vector in $d=2$ space dimensions, and let $t \in \mathbb{R}^+$ be the time coordinate. The Euler equations of hydrodynamics can be written in the updated Lagrangian form using the material derivative $\mde{}{t}$ as follows:
\begin{subequations}
	\label{eqn.PDE_energy}
	\begin{align}
		&\rho \mde{\tau}{t}-\nabla \cdot \vv=0,    \label{eqn.cl1}\\
		&\rho \mde{\vv}{t}+\nabla p=\mathbf{0},    \label{eqn.cl2}\\
		&\rho \mde{\ee}{t}+\nabla \cdot (p \vv)=0 \label{eqn.cl3}.
	\end{align}
\end{subequations}
Here, $\rho$ is the fluid density and $\tau=1/\rho$ is the specific volume, $\vv \in \mathbb{R}^d$ is the velocity vector, $p$ denotes the fluid pressure and $\ee=\varepsilon+\halb \vv^2$ is the specific total energy. Furthermore, in the Lagrangian framework, the system is supplemented with the trajectory equation for the fluid particles, that is
\begin{equation}
	\mde{\x}{t} = \vv, \qquad \x(0)=\mathbf{X},
	\label{eqn.trajODE}
\end{equation}
where $\mathbf{X}$ is the Lagrangian coordinate defined at time $t=0$ corresponding to the Eulerian coordinate $\x$ for $t>0$. The system is closed by an equation of state (EOS) that provides the specific internal energy $\varepsilon$ in terms of the specific volume $\tau$ and specific entropy $\S$, hence $\varepsilon = \varepsilon(\tau,\S)$. In addition, $\varepsilon(\tau,\S)$ being convex implies that $\S(\tau,\varepsilon)$ is concave, thus dealing with the physical entropy. To enforce thermodynamic stability \cite{Menikoff198975,Godlewski}, the specific internal energy is assumed to be convex with respect to $\tau$ and $\S$. Consequently, if the specific internal energy is a thermodynamic potential, the pressure $p$ and the temperature $\temp$ can be determined by the complete equation of state
\begin{equation}
	p(\tau,\S) = -\left( \frac{\partial \varepsilon}{\partial \tau}\right)_{\S}, \qquad \temp(\tau,\S) = \left( \frac{\partial \varepsilon}{\partial \S}\right)_{\tau},
	\label{eqn.en_pot}
\end{equation} 
where the absolute temperature is assumed to be strictly positive ($\temp>0$). Thanks to the convexity of the EOS we define the isentropic sound speed $a$ as
\begin{equation}
	\frac{a^2}{\tau^2} = -\frac{\partial p}{\partial \tau} = \frac{\partial^2 \varepsilon}{\partial \tau^2} > 0.
\end{equation}
The ideal gas EOS fulfills these assumptions, and it simply writes
\begin{equation}
	\varepsilon = \frac{e^{\S/c_v}}{\tau^{\gamma-1} \, (\gamma-1)},
	\label{eqn.eos}
\end{equation}
with $\gamma=c_p/c_v$ being the polytropic index of the gas given as the ratio between the specific heat at constant pressure and volume, respectively. The Lagrangian hydrodynamics equations \eqref{eqn.PDE_energy} have the following set of real eigenvalues:
\begin{equation}
	\boldsymbol{\lambda} = (-a,\mathbf{0},a), \qquad a=\sqrt{\gamma \, p \, \tau},
	\label{eqn.eigenval}
\end{equation}   
with $a$ being the isentropic speed of sound.
The convexity of $\varepsilon$ implies that the specific total energy $\ee$ is also convex. The Gibbs relation stems from the EOS:
\begin{eqnarray}
	\text{d}\ee &=& \frac{\partial \ee}{\partial \tau} \, \text{d}\tau + \frac{\partial \ee}{\partial \vv} \, \text{d}\vv + \frac{\partial \ee}{\partial \S} \, \text{d}\S \nonumber \\
	&=&  -p \, \text{d}\tau + \vv \, \text{d}\vv + \temp \, \text{d}\S = \w \cdot \text{d}\q, \qquad \w= \frac{\partial \ee}{\partial \q} = (-p,\vv,\temp)^\top.
	\label{eqn.de}
\end{eqnarray}
where $\q=(\tau,\vv,\ee)^\top$ denotes the vector of conserved variables, while $\w$ is the vector of dual or Godunov variables \cite{God1961} of the energy potential. Similarly, from the Gibbs relation \eqref{eqn.de} one has
\begin{equation}
	\text{d}\S = \frac{1}{\temp} \left( p \, \text{d}\tau - \vv \, \text{d}\vv + \text{d}\ee \right) = \w^* \cdot \text{d}\q, \qquad \w^* = \frac{\partial \S}{\partial \q} = \frac{1}{\theta}(p,-\vv,1)^{\top},
	\label{eqn.dS}
\end{equation}
with $\w^{*}$ being the vector of dual variables of the entropy potential. Assuming that the flow variables are smooth, the governing equations \eqref{eqn.PDE_energy} admit one additional conservation law for entropy, that is derived by dot-multiplying the Lagrangian hydrodynamics system with the dual variables $\w^{*}$ defined by \eqref{eqn.dS}:
\begin{equation}
	\rho \mde{\S}{t} + \frac{1}{\temp} \left( -p \, \nabla \cdot \vv - \vv \cdot \nabla p + \nabla \cdot (p\vv)  \right) = 0.
\end{equation}
Therefore, entropy is conserved under the condition
\begin{equation}
	-p \, \nabla \cdot \vv - \vv \cdot \nabla p + \nabla \cdot (p\vv) = 0,
	\label{eqn.prod_rule}
\end{equation}
which is trivial if smooth flows are considered because this is the algebraic identity of the product rule. In case of discontinuous flows, entropy production takes place in order to discriminate the physical admissible solution characterized by an entropy increase, hence leading to an entropy inequality. Differently from most of the existing solvers for Lagrangian hydrodynamics \cite{Maire2007,Maire2009,maire_loubere_vachal10,Despres2005,Despres2009,ShashkovCellCentered,Burton2015}, we aim at directly solving the entropy inequality instead of the energy equation, hence obtaining the energy equation as a consequence of the discrete form of the following conservation laws:
\begin{subequations}
	\label{eqn.PDE_entropy}
	\begin{align}
		&\rho \mde{\tau}{t}-\nabla \cdot \vv=0,    \label{eqn.cl_tau}\\
		&\rho \mde{\vv}{t}+\nabla p=\mathbf{0},    \label{eqn.cl_v}\\
		&\rho \mde{\S}{t}= \Pi \, \geq 0, \label{eqn.cl_S}
	\end{align}
\end{subequations}
with $\Pi$ being a suitable entropy production term that must be non-negative to obtain physically admissible solutions. By dot-multiplying the above equations \eqref{eqn.cl_tau}-\eqref{eqn.cl_S} with the dual variables $\w$ defined in \eqref{eqn.de}, we retrieve the energy equation
\begin{eqnarray}
	\rho \mde{\ee}{t} + p \nabla \cdot \vv  + \vv \cdot \nabla p &=& \temp \Pi, \nonumber \\
	\rho \mde{\ee}{t} + \nabla \cdot (p \vv) &=& \temp \left( \Pi - \frac{1}{\theta} \left( p \, \nabla \cdot \vv + \vv \cdot \nabla p - \nabla \cdot (p\vv) \right)\right),
	\label{eqn.en_consequence}
\end{eqnarray}
which exactly corresponds to \eqref{eqn.cl3} if $\Pi=\left( p \, \nabla \cdot \vv + \vv \cdot \nabla p - \nabla \cdot (p\vv) \right)/\theta$, thus for smooth solutions we retrieve the entropy preserving condition \eqref{eqn.prod_rule}. Therefore, the relation \eqref{eqn.prod_rule} also ensures that energy is conserved under the assumption of smooth solutions. In the next section, we design a numerical method that can preserve this thermodynamic compatibility also at the discrete level and also in the presence of discontinuities.

\section{Numerical method} \label{sec.method}
The Euclidean space $\Omega(t)$ with boundary $\partial \Omega(t)$ in $d=2$ space dimensions is discretized with a set of non-overlapping triangular control volumes $\omega_c(t)$ of volume $|\omega_c(t)|$. 
\begin{figure}[!htbp]
	\begin{center}
		\begin{tabular}{lll}        
			\includegraphics[trim=10 0 150 10,clip,width=0.30\textwidth]{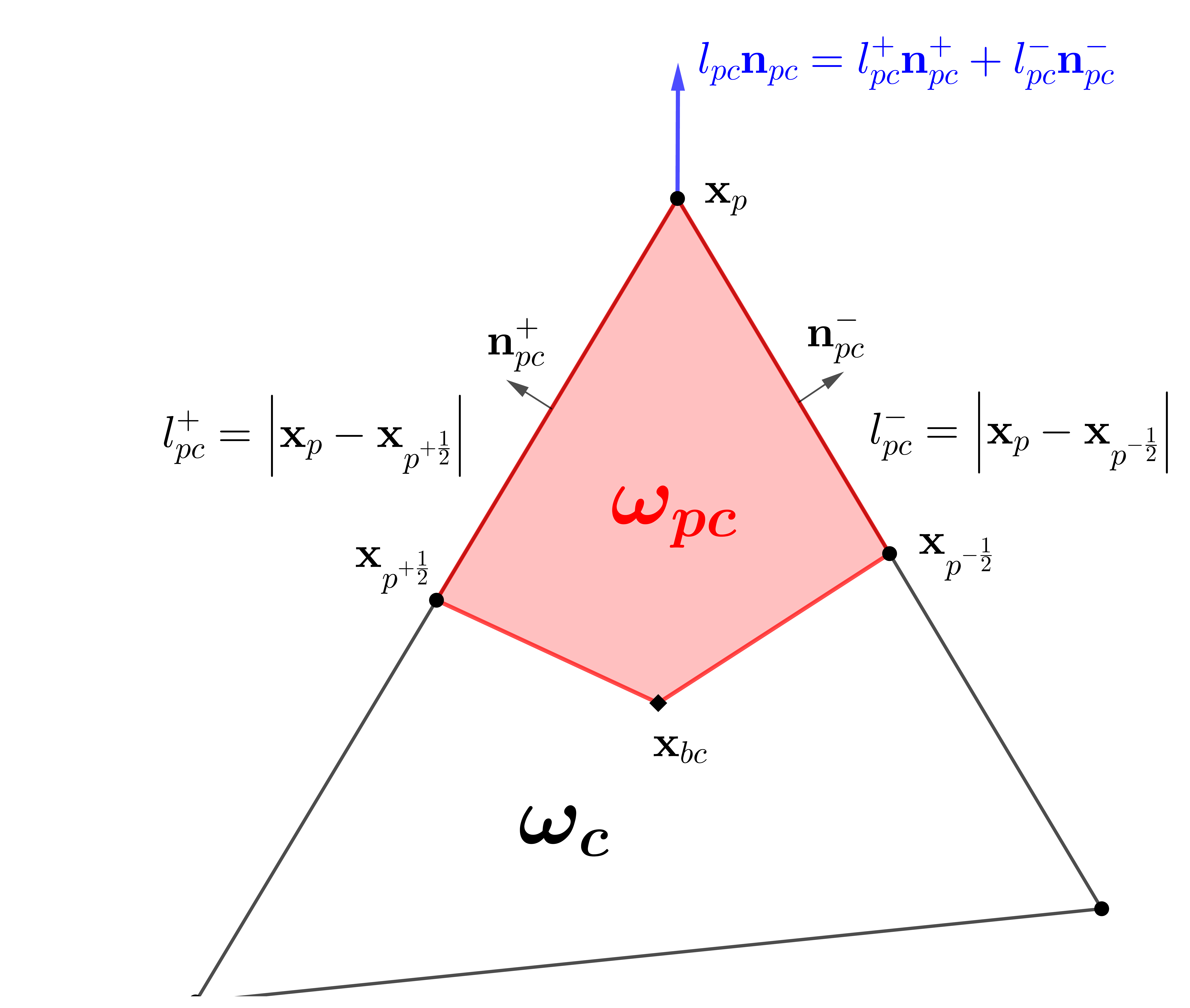} &
			\includegraphics[trim=100 100 100 100,clip,width=0.30\textwidth]{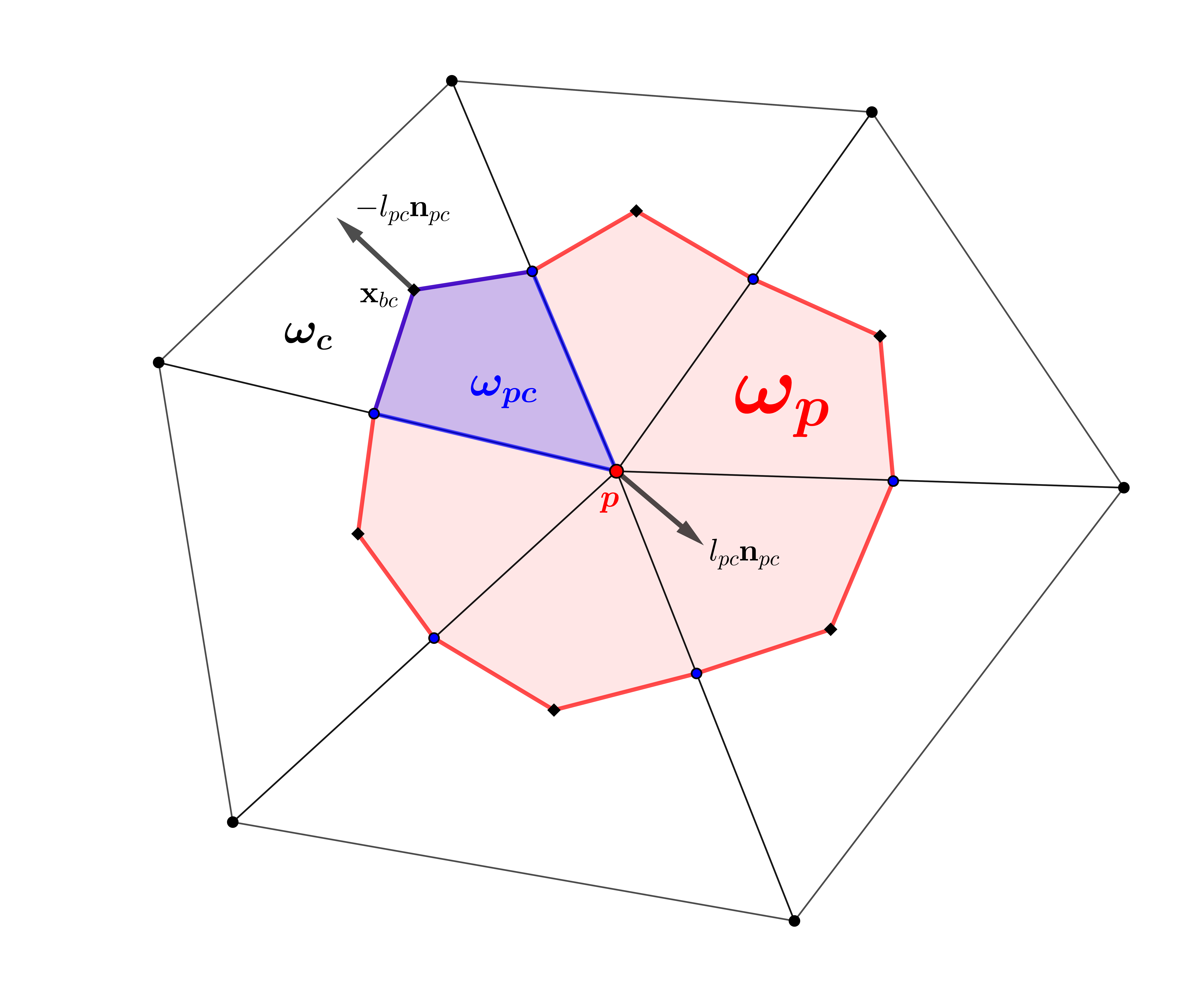} &
			\includegraphics[trim=80 10 30 10,clip,width=0.33\textwidth]{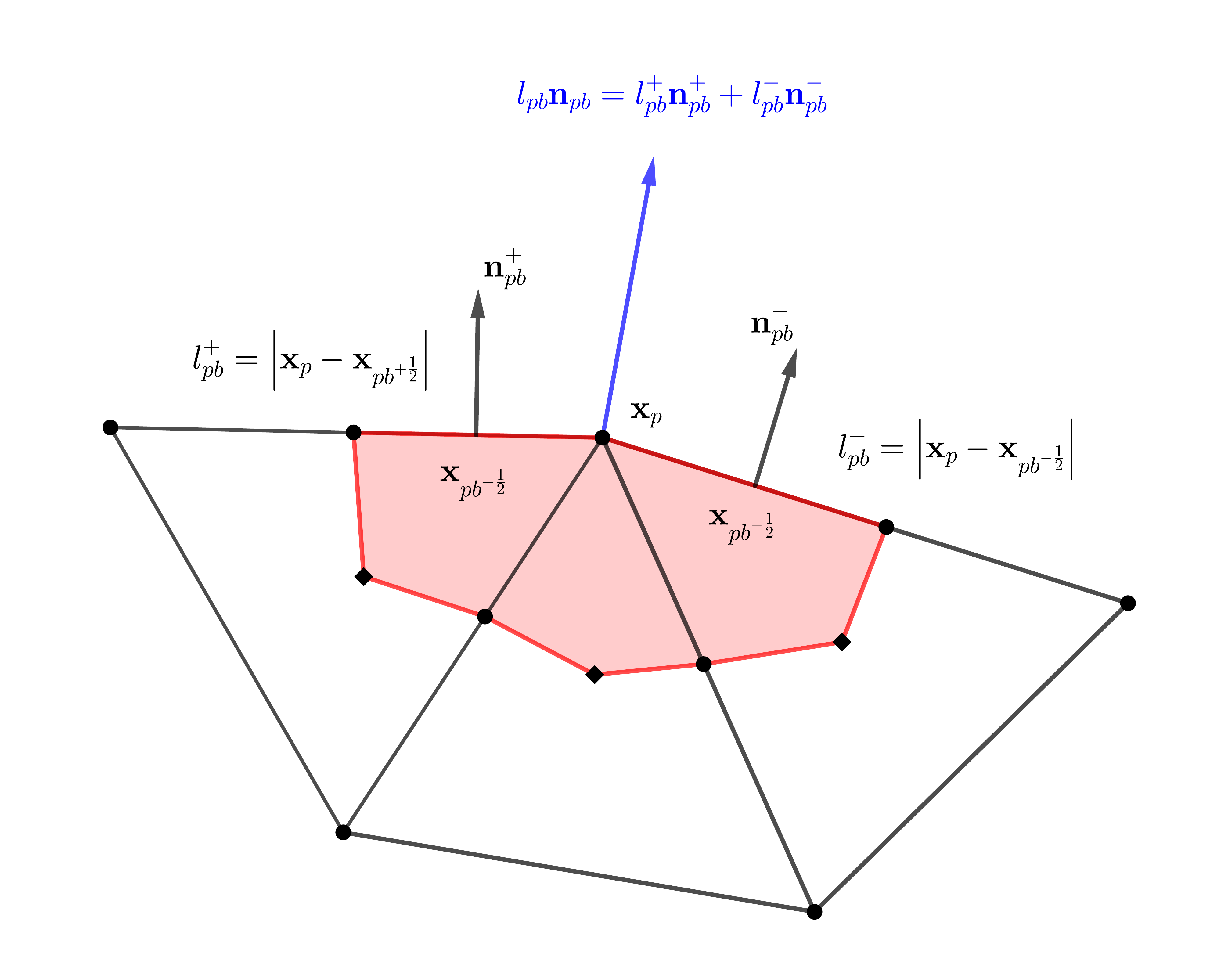}\\
		\end{tabular} 
		\caption{Notation for the cell $\omega_c$ (left), definition of the dual cell $\omega_p$ (middle) and notation for the boundary cell located on $\partial \Omega$.}
		\label{fig.notation}
	\end{center}
\end{figure}
Notice that the geometry related quantities are time dependent due to the mesh motion, while the topology of the grid remains unchanged. A generic vertex of the mesh is denoted by $p$, while the sub-cell $\omega_{pc}$ refers to the portion of cell $\omega_c$ attached to one of its vertex $p$, as depicted in Figure \ref{fig.notation}. The set of vertexes belonging to cell $\omega_c$ is referred to with $\mathcal{P}(c)$, while the set of cells sharing node $p$ is indicated with $\mathcal{C}(p)$. The sub-cell $\omega_{pc}$ is defined by connecting the vertex with position $\x_p$, the cell barycenter of coordinates $\x_{bc}=1/(d+1) \, \sum_{p \in \mathcal{P}(c)}\x_p$, and the left and right midpoints of the edges impinging on node $p$, namely $p^{-\halb}$ and $p^{+\halb}$, respectively. The half lengths of these edges are given by $l_{pc}^-=\left| \x_p - \x_p^{-\halb} \right|$ and $l_{pc}^+=\left| \x_p - \x_p^{+\halb} \right|$, and the corner normal is computed as
\begin{equation}
	\lnpc = l_{pc}^+ \n_{pc}^+ + l_{pc}^- \n_{pc}^-.
	\label{eqn.lpcnpc}
\end{equation}
By construction the corner vectors satisfy the fundamental geometrical identity
\begin{equation}
	\sum \limits_{p \in \mathcal{P}(c)} \lnpc = \mathbf{0},
	\label{eqn.gauss}
\end{equation}
which is a consequence of the Gauss theorem. At the aid of Figure \ref{fig.notation}, we introduce the dual cell $\omega_p$ given by the union of all sub-cells sharing a generic node $p$:
\begin{equation}
	\omega_p = \bigcup \limits_{c \in \mathcal{C}(p)} \omega_{pc}.
	\label{eqn.dualcell}
\end{equation}
Notice that the outward pointing corner vectors of $\omega_p$ are simply the cell corner vectors $\lnpc$ with opposite sign for all $c \in \mathcal{C}(p)$, thus the following identity holds
\begin{equation}
	\sum \limits_{c \in \mathcal{C}(p)} -\lnpc = \mathbf{0}.
	\label{eqn.gaussd}
\end{equation}
The time coordinate $t$ is defined in the interval $[0;t_f]$ and is approximated by a sequence of discrete points $t^n$ such that
\begin{equation}
	t^{n+1}=t^n + \dt, 
\end{equation}
where the time step $\dt$ is computed according to a CFL stability condition as
\begin{equation}
	\dt \leq \text{CFL} \, \min \limits_{c} \frac{|\omega_c|}{a_c}.
	\label{eqn.timestep}
\end{equation}

\subsection{Semi-discrete entropy conservative Lagrangian schemes (ECL)} \label{ssec.EC}
We first design a cell-centered finite volume scheme 
that is exactly entropy conservative, thus we set $\Pi=0$ in the entropy inequality \eqref{eqn.cl_S}. Let $m_c=|\omega_c(t)|/\tau_c(t)$ be the mass of the cell, which remains constant in the Lagrangian framework, and let the mass averaged value of a generic quantity $\phi(\x,t)$ be defined as
\begin{equation}
	\phi_c = \frac{1}{m_c} \, \int \limits_{\omega_c(t)} \rho \phi \, \text{d}\x.
\end{equation}
By integrating system \eqref{eqn.PDE_entropy} over the cell $\omega_c(t)$ and using the Reynolds transport formula, the semi-discrete finite volume scheme in fluctuation form reads 
\begin{subequations}
	\label{eqn.SD_entropy}
	\begin{align}
		&m_c \mde{\tau_c}{t}- \sum \limits_{p \in \mathcal{P}(c)} \lnpc \cdot (\vv_p-\vv_c)=0,    \label{eqn.sd_tau}\\
		&m_c \mde{\vv_c}{t}+\sum \limits_{p \in \mathcal{P}(c)} \lnpc \, (\tilde{p}_p-p_c)=\mathbf{0},    \label{eqn.sd_v}\\
		&m_c \mde{\S_c}{t}= 0. \label{eqn.sd_S}
	\end{align}
\end{subequations}    
%
The motion of the computational grid is governed by the trajectory equation \eqref{eqn.trajODE} that is discretized at each node of the mesh as 
\begin{equation}
	\mde{\x_p}{t} = \vv_p, \qquad \x_p(0)= \mathbf{X}(0).
	\label{eqn.xp}
\end{equation}
To rigorously ensure the geometric conservation law (GCL), the flux approximation of the continuity equation \eqref{eqn.sd_tau} must provide a time rate of change of the cell volume that is fully compatible with the mesh displacement given by \eqref{eqn.xp}. This comes from the fact that $\partial |\omega_c|/\partial \x_p=\lnpc$ and thus \eqref{eqn.sd_tau} corresponds to $\text{d}|\omega_c|/\text{d}t$ \cite{CCLAM_hyper}. In other words, the volume of the cell given by $\tau_c^{-1}$ must be equal to the volume of the cell computed from the vertex coordinates. Therefore, the node velocity in \eqref{eqn.sd_tau} and \eqref{eqn.xp} must be the same, and one is not allowed to modify this requirement, which is the compliance with the GCL. However, we can correct the nodal pressure flux $\tilde{p}_p$ to be the sought thermodynamically compatible numerical flux. Following the general framework of entropy conservative schemes proposed by Abgrall in \cite{Abgrall2018}, we assume that the entropy conservative numerical flux writes
\begin{equation}
	\label{eqn.thc_flux}
	\lnpc  \, \tilde{p}_p   = \lnpc \,  p_p + \| \lnpc \| \, \alpha_p \, (\vv_c- \vv_p),
\end{equation}  
with $p_p$ an averaged node pressure that is not necessarily thermodynamically compatible and $\alpha_p$ a scalar nodal correction factor that multiplies the jump in the thermodynamic dual variables and that eventually corrects the flux $p_p$ to ensure thermodynamic consistency. 
Therefore, the semi-discrete entropy conservative scheme is given by
\begin{subequations}
	\label{eqn.SD_thc}
	\begin{align}
		&m_c \mde{\tau_c}{t} - \sum \limits_{p \in \mathcal{P}(c)} \lnpc \cdot ({\vv}_p-\vv_c) = 0,    \label{eqn.sdthc_tau}\\
		&m_c \mde{\vv_c}{t} + \sum \limits_{p \in \mathcal{P}(c)} \lnpc \, ({p}_p-p_c) + \| \lnpc \| \, \alpha_p \, (\vv_c- \vv_p) = \mathbf{0},    \label{eqn.sdthc_v}\\
		&m_c \mde{\S_c}{t} = 0. \label{eqn.sdthc_S}
	\end{align}
\end{subequations} 

First, we need to determine the non-compatible fluxes $\vv_p$ and $p_p$ invoking conservation principles. A consistent condition implies that the sum of the fluctuations around a node must be equal to the sum of the fluxes across the faces defining the dual cell $\omega_p$, hence we require
\begin{equation}
	\label{eqn.conservation}
	\sum \limits_{c \in \mathcal{C}(p)} \lnpc \, (p_p-p_c) + \alpha_p  \sum \limits_{c \in \mathcal{C}(p)} \| \lnpc \| (\vv_c-\vv_p) = - \!\! \sum \limits_{c \in \mathcal{C}(p)} \lnpc \, p_c \,  = \!\! \int \limits_{\partial \omega_p(t)} \!\! p \cdot \n \, \text{ds},
\end{equation}
where the change of sign is due to the opposite corner vectors defined on the dual cell according to \eqref{eqn.gaussd}. By canceling the terms appearing on both sides of \eqref{eqn.conservation} and recalling the discrete Gauss theorem on the dual cell \eqref{eqn.gaussd}, we obtain the nodal velocity which satisfies conservation. Although the nodal pressure $p_p$ is in principle without any constraints, we use the same definition of the conservative velocity flux, thus we obtain
\begin{equation}
	\vv_p = \frac{\sum \limits_{c \in \mathcal{C}(p)} \| \lnpc \| \, \vv_c}{\sum \limits_{c \in \mathcal{C}(p)} \| \lnpc \|}, \qquad p_p = \frac{\sum \limits_{c \in \mathcal{C}(p)} \| \lnpc \| \, p_c}{\sum \limits_{c \in \mathcal{C}(p)} \| \lnpc \|}.
	\label{eqn.pp_vp}
\end{equation}

Next, we need to provide thermodynamic compatibility by imposing the condition \eqref{eqn.prod_rule} around a node $p$, thus we compute the dot product of the discrete dual variables \eqref{eqn.de} with the fluxes of the semi-discrete scheme \eqref{eqn.SD_thc}: 
\begin{eqnarray}
	&&  \qquad  \sum \limits_{c \in \mathcal{C}(p)} \! \lnpc \cdot \left( p_c \, (\vv_p-\vv_c) + \vv_c \, (p_p-p_c) \right) 
+ \alpha_p \!\! \sum \limits_{c \in \mathcal{C}(p)} \| \lnpc \| \, \vv_c \cdot (\vv_c-\vv_p)  \nonumber \\ 
	&& \qquad = - \sum \limits_{c \in \mathcal{C}(p)} \! \lnpc \cdot p_c \, \vv_c  = \int \limits_{\partial \omega_p(t)} p \, \vv \cdot \n \, \text{ds}.   
	\label{eqn.Ethc}
\end{eqnarray}
The term multiplied by the scalar factor $\alpha_p$ can be conveniently rearranged in quadratic form as follows:
\begin{equation}
	\label{eqn.alphaterms2}	
	  \sum \limits_{c \in \mathcal{C}(p)} \| \lnpc \| \, (\vv_c - \vv_p)^2 = \sum \limits_{c \in \mathcal{C}(p)} \| \lnpc \| \, \vv_c \, (\vv_c-\vv_p) - \vv_p \, \underbrace{\sum \limits_{c \in \mathcal{C}(p)} \| \lnpc \| \, (\vv_c-\vv_p)}_{=0},
\end{equation}
where the last term on the right hand side is zero due to the definition of the nodal fluxes \eqref{eqn.pp_vp}. Inserting the above expression \eqref{eqn.alphaterms2} into the thermodynamically compatible condition \eqref{eqn.Ethc}, allows the scalar correction factor $\alpha_p$ to be computed:
\begin{equation}
	\alpha_p = \frac{\sum \limits_{c \in \mathcal{C}(p)} \!\!\! \lnpc \cdot (p_c \, \vv_c - \vv_p \, p_c - p_p \, \vv_c)}{\sum \limits_{c \in \mathcal{C}(p)} \!\! \| \lnpc \| \, (\vv_c - \vv_p)^2 } = \frac{|\omega_p| \left( -\nabla_p \! \cdot \! (p\vv) + \vv_p \! \cdot \! \nabla_p p + p_p \nabla_p \! \cdot \! \vv \right)}{\sum \limits_{c \in \mathcal{C}(p)} \| \lnpc \| \, (\vv_c - \vv_p)^2 }
	\label{eqn.alpha}
\end{equation}
where $\nabla_p \cdot$ and $\nabla_p$ are the discrete divergence and gradient operators defined on the dual cell \cite{Maire2020}. It is interesting to notice that the condition \eqref{eqn.prod_rule} can also be interpreted as a consistency condition on the dual cell (similar to \eqref{eqn.conservation}) imposed on the contraction of the fluctuations with the dual variables $\w$, that indeed must be equal to the energy flux $\nabla \cdot (p\vv)$. In the case of vanishing denominator in \eqref{eqn.alpha}, that is when its numerical value is lower than $10^{-20}$, the correction factor is set to $\alpha_p=0$. Notice that the thermodynamically compatible scheme \eqref{eqn.SD_thc} exhibits no numerical dissipation since the fluxes \eqref{eqn.pp_vp} are given by geometry weighted averages that can be seen as purely central fluxes. Nevertheless, the semi-discrete scheme \eqref{eqn.SD_thc} is nonlinearly stable in the energy norm, as we will prove later.

\begin{theorem} \label{th.NLstability}
Assuming impermeable boundary conditions $\int \limits_{\partial \Omega} \vv \cdot \n \, \textnormal{ds}= 0$,
the semi-discrete scheme \eqref{eqn.SD_thc} with nodal fluxes defined by \eqref{eqn.pp_vp} and the correction factor given by \eqref{eqn.alpha} is nonlinearly stable in the energy norm in the sense that we have
\begin{equation}
	\int \limits_{\Omega} \mde{\ee}{t} \, \text{d}\x = 0.
	\label{eqn.NLstability}
\end{equation}	
\end{theorem}

\begin{proof}
Mimicking the procedure at the continuous level yielding \eqref{eqn.en_consequence}, let us compute the semi-discrete energy conservation law by contracting the semi-discrete scheme \eqref{eqn.SD_thc} with the dual variables $\w$ given by \eqref{eqn.de}: 
\begin{eqnarray}
m_c \mde{\ee_c}{t} &=& \sum \limits_{p \in \mathcal{P}(c)} \lnpc \cdot \left( p_c \, (\vv_c-\vv_p) + \vv_c \, (p_c-p_p) \right) \nonumber \\
&-& \sum \limits_{p \in \mathcal{P}(c)} \| \lnpc \| \, \alpha_p \, \vv_c \cdot (\vv_c-\vv_p). 	
\end{eqnarray} 
The above equation can now be summed over all the cells which discretize the computational domain, yielding
\begin{eqnarray}
	\sum \limits_{c} \, m_c \mde{\ee}{t} &=& \sum \limits_{c}  \, \sum \limits_{p \in \mathcal{P}(c)} \lnpc \cdot \left( p_c \, (\vv_c-\vv_p) + \vv_c \, (p_c-p_p) \right) \nonumber \\
	&-& \sum \limits_{c} \, \sum \limits_{p \in \mathcal{P}(c)} \| \lnpc \| \, \alpha_p \, \vv_c \cdot  (\vv_c-\vv_p) \nonumber \\
	&=& \sum \limits_{p}  \, \sum \limits_{c \in \mathcal{C}(p)} \lnpc \cdot \left( p_c \, (\vv_c-\vv_p) + \vv_c \, (p_c-p_p) \right) \nonumber \\
	&-& \sum \limits_{p} \, \sum \limits_{c \in \mathcal{C}(p)} \| \lnpc \| \, \alpha_p \, \vv_c \cdot  (\vv_c-\vv_p) = 0, 	
\end{eqnarray} 
where we have switched the summation over all the cells and all the nodes of each cell with the summation over all the nodes and all the cells surrounding each node. Notice that the right hand side of the above equation is exactly the right hand side of the thermodynamic compatibility condition \eqref{eqn.Ethc}, which together with the definition of $\alpha_p$ \eqref{eqn.alpha} leads to the last equality, 
hence retrieving the sought nonlinear stability in the energy norm \eqref{eqn.NLstability}.
\end{proof}		

\subsection{Boundary conditions} \label{ssec.BCs}
In the Lagrangian framework one has to deal with two types of boundary conditions, imposing either the pressure or the normal component of the velocity. 
Let $p$ be a node lying on the boundary of the domain $\partial \Omega$, as depicted in Figure \ref{fig.notation}. The corner vector on the boundary is given by 
\begin{equation}
	\lnpb = l_{pb}^{+} \n_{pb}^{+} + l_{pb}^{-} \n_{pb}^{-},
	\label{eqn.lpcnpcb}
\end{equation}
where the outward normals $l_{pb}^{+} \n_{pb}^{+}$ and $l_{pb}^{-} \n_{pb}^{-}$ are referred to the boundary sides impinging on node $p$. These normals are essentials to satisfy the discrete Gauss theorem also for a boundary node, and they are indeed linked to the corner vectors with the following relation:
\begin{equation}
	\sum \limits_{c \in \mathcal{C}(p)} \lnpc = \lnpb.
\label{eqn.gaussb}
\end{equation}
Particular care must be devoted to design the boundary conditions in such a way that the GCL property of the scheme is not spoiled. Here, we consider pressure and wall boundaries.

\paragraph{Pressure boundary condition} For a boundary node the consistent condition \eqref{eqn.conservation} would write
\begin{equation}
	\label{eqn.conservation_bnd_prs}
	-\sum \limits_{c \in \mathcal{C}(p)} \lnpc \, (p_c-p_p) + \alpha_p  \sum \limits_{c \in \mathcal{C}(p)} \| \lnpc \| (\vv_c-\vv_p) = -\sum \limits_{c \in \mathcal{C}(p)} \lnpc \, p_c + \lnpb \, p_b,
\end{equation}
where $p_b$ is the prescribed nodal pressure. We recall that the definition of $p_p$ is not constrained under any conservation requirement, thus the choice made in \eqref{eqn.pp_vp} can be modified. In order to maintain the same definition of $\vv_p$ in \eqref{eqn.pp_vp}, which ensures conservation, we impose the following boundary pressure:
\begin{equation}
	p_p = p_b.
	\label{eqn.pb}
\end{equation}
Inserting the above definition into \eqref{eqn.conservation_bnd_prs} and using the geometry identity \eqref{eqn.gaussb}, we retrieve the same conservation condition on the velocity given by \eqref{eqn.pp_vp}. The thermodynamic compatibility condition \eqref{eqn.Ethc} and thus the computation of the correction factor $\alpha_p$ \eqref{eqn.alpha} remain unchanged.

\paragraph{Wall boundary condition} In order to prescribe the zero normal velocity at a boundary node, let us decompose the velocity vector of a boundary cell into the normal and tangential contributions with respect to the boundary corner vector \eqref{eqn.lpcnpcb}:
\begin{equation}
 \vv_c := \vv_{c}^\perp + \vv_{c}^\parallel = (\vv_c \cdot \lnpb) \, \lnpb + \left( \vv_c - \frac{\vv_c \cdot \lnpb}{\| \lnpb \|}\, \lnpb \right).  	
\end{equation}
Now, since the normal contribution is prescribed as boundary condition, we only have the freedom to determine the tangential velocity vector, relying on the consistency condition \eqref{eqn.conservation}, thus we obtain
\begin{eqnarray}
	\label{eqn.conservation_bnd_vel}
	-\sum \limits_{c \in \mathcal{C}(p)} \lnpc \, (p_c-p_p) &+& \alpha_p  \sum \limits_{c \in \mathcal{C}(p)} \| \lnpc \| \left( (\vv_{c}^\perp + \vv_{c}^\parallel) -\vv_p \right) \nonumber \\
	&=& -\sum \limits_{c \in \mathcal{C}(p)} \lnpc \, p_c + \lnpb \, p_b.
\end{eqnarray} 
The boundary nodal pressure can be arbitrarily chosen because we are imposing a velocity boundary condition, therefore we simply set $p_b=p_p$. For wall boundaries obviously $\vv_{c}^\perp=0$ holds, hence the velocity vector that is compatible with the prescribed boundary condition and satisfies conservation is eventually given by
\begin{equation}
\vv_p = \frac{\sum \limits_{c \in \mathcal{C}(p)} \| \lnpc \| \, \vv_{c}^\parallel}{\sum \limits_{c \in \mathcal{C}(p)} \| \lnpc \|}.
\end{equation}
This new definition of the node velocity must be taken into account also in the thermodynamically compatible correction \eqref{eqn.Ethc}, that now becomes  
\begin{eqnarray}
	\sum \limits_{c \in \mathcal{C}(p)} \lnpc \cdot p_c \, \vv_c &=& \phantom{\alpha_p \, } \sum \limits_{c \in \mathcal{C}(p)} \lnpc \cdot \left( p_c \, (\vv_c-\vv_p) + \vv_c \, (p_c-p_p) \right) \nonumber \\
	&-& \alpha_p \, \sum \limits_{c \in \mathcal{C}(p)} \| \lnpc \| \, (\vv_{c}^\parallel - \vv_p)^2, 
	\label{eqn.Ethc_bnd_vel}
\end{eqnarray}
from which we derive the correction factor $\alpha_p$ that is also consistent with the prescribed wall boundary condition.

\subsection{Semi-discrete entropy stable Lagrangian schemes (ESL)} \label{ssec.ED}
If the flow exhibits discontinuities like shock waves, the entropy inequality \eqref{eqn.cl_S} permits to define the physically admissible solution characterized by an entropy increase. In this case we need to explicitly consider the entropy production term $\Pi$ in our new scheme, which directly solves the entropy inequality. The production term must exactly account for the numerical dissipation introduced in the scheme, hence we need to: i) supplement the scheme with suitable numerical diffusion and, ii) quantify the amount of numerical viscosity to obtain a thermodynamically compatible entropy stable scheme.

The dissipative scheme is based on the EUCCLHYD finite volume method \cite{EUCCLHYD,Maire2009}, which is here rewritten in fluctuation form for the governing equations \eqref{eqn.PDE_entropy}:
\begin{subequations}
	\label{eqn.EUC}
	\begin{align}
		&m_c \mde{\tau_c}{t} - \sum \limits_{p \in \mathcal{P}(c)} \lnpc \cdot ({\vv}_p^*-\vv_c) = 0,    \label{eqn.euc_tau}\\
		&m_c \mde{\vv_c}{t} + \sum \limits_{p \in \mathcal{P}(c)} \lnpc \, ({p}_{pc}^*-p_c) = \mathbf{0},    \label{eqn.euc_v}\\
		&m_c \mde{\S_c}{t} = \sum \limits_{p \in \mathcal{P}(c)} \Pi_{pc} \geq 0. \label{eqn.euc_S}
	\end{align}
\end{subequations}
where $\lnpc \cdot {\vv}_p^*$ and $\lnpc {p}_p^*$ are the consistent numerical flux functions for the continuity and the momentum equations, respectively. In order to satisfy an entropy inequality within the cell $\omega_c$, the expression of the pressure flux, also known as the sub-cell force \cite{Maire2011}, is given by
\begin{equation}
\lnpc {p}_{pc}^* = \lnpc p_c - \Mpc (\vv_p^*-\vv_c),	
\label{eqn.fpc}
\end{equation}	
with the symmetric and positive definite corner viscosity matrix \cite{Maire2011}
\begin{equation}
	\Mpc = \frac{a_c}{\tau_c} \left( l_{pc}^- \, (\n_{pc}^- \otimes \n_{pc}^-) + l_{pc}^+ \, (\n_{pc}^+ \otimes \n_{pc}^+) \right).
	\label{eqn.Mpc}
\end{equation}
The node velocity $\vv_p^*$ is determined such that momentum conservation is ensured. Similar to the entropy conservative scheme, we apply the consistency condition \eqref{eqn.conservation} on the dual cell $\omega_p$ to the semi-discrete scheme \eqref{eqn.EUC}, thus obtaining
\begin{equation}
	\sum \limits_{c \in \mathcal{C}(p)} \lnpc ({p}_{pc}^* - p_c ) = -\sum \limits_{c \in \mathcal{C}(p)} \lnpc p_c = \int \limits_{\partial \omega_p(t)} p\cdot \n \, \textit{ds},
	\label{eqn.fpc_balance}
\end{equation}
which requires that the sum of the sub-cell forces around a node vanishes. It is worth noticing that exactly the same relation is obtained in \cite{EUCCLHYD} by invoking conservation principles of energy and momentum over the whole computational domain. Inserting the definition of the sub-cell force \eqref{eqn.fpc} into the above equation yields a linear system for the nodal velocity:
\begin{equation}
	\Mp \vv_p^* = \sum \limits_{c \in \mathcal{C}(p)} \lnpc p_c + \Mpc \vv_c, \qquad \Mp=\sum \limits_{c \in \mathcal{C}(p)} \Mpc,
	\label{eqn.NS}
\end{equation}
which is typically referred to as nodal solver. To comply with the GCL, the nodes of the grid are moved according to \eqref{eqn.trajODE} using the velocity $\vv_p^*$.
In order to quantify the production term in \eqref{eqn.euc_S}, we dot-multiply  the semi-discrete finite volume scheme \eqref{eqn.EUC} by the vector of dual variables $\w$ defined in \eqref{eqn.de}, hence obtaining
\begin{equation}
	\Pi_{pc} = \frac{1}{\temp_c} \, (\vv_c-\vv_p^*) \, \Mpc \, (\vv_c-\vv_p^*) \geq 0, 
	\label{eqn.prod}
\end{equation} 
which is the non-negative entropy production term needed in eqn. \eqref{eqn.euc_S} of the scheme. It is non-negative since the temperature $\theta_c \geq 0$ is assumed to be non-negative and because the matrix $\Mpc$ is symmetric positive definite. Hence, our new scheme satisfies the entropy inequality
\textit{by construction}, which is a unique feature of the novel cell-centered Lagrangian scheme introduced in this paper. 

\begin{theorem} \label{th.NLstability_ES}
	Assuming impermeable wall boundary conditions $\int \limits_{\partial \Omega} \vv \cdot \n \, \textnormal{ds}= 0$,
	the semi-discrete scheme \eqref{eqn.EUC} with nodal fluxes defined by \eqref{eqn.fpc} and \eqref{eqn.NS} is nonlinearly stable in the energy norm in the sense that we have
	\begin{equation}
		\int \limits_{\Omega} \mde{\ee}{t} \, \text{d}\x = 0.
		\label{eqn.NLstability_EC}
	\end{equation}	
\end{theorem}

\begin{proof}
	Similar to the proof of Theorem \ref{th.NLstability}, we proceed by computing the semi-discrete energy conservation law as the dot-product of the semi-discrete scheme \eqref{eqn.EUC} with the dual variables $\w$ given by \eqref{eqn.de}:
	\begin{eqnarray}
		m_c \mde{\ee_c}{t} &=& -\sum \limits_{p \in \mathcal{P}(c)} p_c \,  \lnpc \cdot  (\vv_p^*-\vv_c)
		                       + \sum \limits_{p \in \mathcal{P}(c)} (\vv_c-\vv_p^*) \, \Mpc \, (\vv_c-\vv_p^*)  \nonumber \\
		&\phantom{=}& - \sum \limits_{p \in \mathcal{P}(c)} \vv_c \cdot \lnpc p_c - \vv_c \cdot \Mpc (\vv_p^*-\vv_c) - \vv_c \cdot \lnpc p_c.
	\end{eqnarray} 
	After some algebra and using the geometrical identity \eqref{eqn.gauss}, the above equation simplifies to
	\begin{eqnarray}
		m_c \mde{\ee_c}{t} &=& - \sum \limits_{p \in \mathcal{P}(c)} p_c \,  \lnpc \cdot  \vv_p^* -\sum \limits_{p \in \mathcal{P}(c)} \vv_p^*  \Mpc \, (\vv_c-\vv_p^*).
		\label{eqn.energy_ES}
	\end{eqnarray}
    Let us now sum the local energy equation \eqref{eqn.energy_ES} over all the cells paving the computational domain, hence obtaining
    \begin{eqnarray}
    	\sum \limits_{c} \, m_c \mde{\ee_c}{t} &=& \sum \limits_{c} \, \sum \limits_{p \in \mathcal{P}(c)} -p_c \,  \lnpc \cdot  \vv_p^* - \vv_p^*  \Mpc \, (\vv_c-\vv_p^*) \nonumber \\
    	&=& -\sum \limits_{p} \, \vv_p^* \, \sum \limits_{c \in \mathcal{C}(p)} p_c \,  \lnpc - \Mpc \, (\vv_p^*-\vv_c) = 0,
    	\label{eqn.energy_ES2}
    \end{eqnarray}
    where we have switched the summation between cells and nodes since we are dealing with finite summations. Notice that the term on the right hand side of the energy balance \eqref{eqn.energy_ES2} is the sub-cell force definition \eqref{eqn.fpc}, whose summation over all the cells surrounding a node is exactly zero because of the conservation property \eqref{eqn.fpc_balance} which the semi-discrete scheme \eqref{eqn.EUC} relies on. Therefore, we retrieve the nonlinear stability \eqref{eqn.NLstability_EC} at the discrete level.
\end{proof}	

\subsection{Adaptive semi-discrete ECL/ESL schemes} \label{ssec.Hybridization}
Because of the non-linearity of the governing equations, the flow can dynamically face situations of discontinuities and/or smooth solutions, that must be properly handled by an adaptive switching in space and time between the entropy conservative and entropy stable scheme. Let $\mathcal{F}(\q)=(\nabla \cdot \vv, -\nabla p, \Pi)^\top$ be the spatial fluxes of the hydrodynamics system \eqref{eqn.PDE_entropy}, and let $\mathcal{F}^{C}_h$ and $\mathcal{F}^{S}_h$ represent its entropy conservative and entropy stable spatial discretization according to \eqref{eqn.SD_thc} and \eqref{eqn.EUC}, respectively. The semi-discrete schemes \eqref{eqn.SD_thc} and \eqref{eqn.EUC} are then hybridized at each time step with a convex combination of the form
\begin{equation}
	m_c \mde{\q_c}{t} = (1-\beta_p) \cdot \mathcal{F}^{C}_h + \beta_p  \cdot \mathcal{F}^{S}_h,
	\label{eqn.HybridScheme}
\end{equation}
where $\beta_p$ is the blending factor. The design of the blending factor is out of the scope of this work, therefore we simply set
\begin{equation}
	\beta_p = \left \{ \begin{array}{ll}
		1 & \textnormal{for discontinuous flows} \\
		0 & \textnormal{for smooth flows}
	\end{array} \right. .
\end{equation} 
To detect discontinuities in the solution, we use either the \textit{a priori} technique designed in \cite{BalsaraFlattener}, which is based on the measure of the discrete divergence of the velocity field, or the \textit{a posteriori} MOOD strategy \cite{MOOD,Boscheri_hyperelast_22,LAM2018}, that identifies the troubled cells relying on a set of physical and numerical admissibility criteria.

\begin{theorem} 
	Assuming a convex equations of state $\varepsilon=\varepsilon(\tau,\S)$, well prepared initial data in the sense that $\tau_c(t=0)>0$ and $S_c(t=0) \geq 0$ and assuming positive cell volumes $|\omega_c(t)|>0$ for all times $t\geq 0$, the semi-discrete scheme \eqref{eqn.HybridScheme} is positivity preserving for pressure, density and temperature, that is
	\begin{equation}
		\rho_c(t) > 0, \qquad p_c(t) \geq 0, \qquad \theta_c(t) \geq 0.
	\end{equation} 	
\end{theorem}	

\begin{proof}
The adaptive scheme \eqref{eqn.HybridScheme} is a convex combination of the semi-discrete schemes \eqref{eqn.SD_thc} and \eqref{eqn.EUC}. Since the Lagrangian element mass $m_c = |\omega_c(0)| / \tau_c(0) > 0$ is strictly positive and constant in time and since the volumes are assumed to be positive, $|\omega_c(t)|>0$, the density $\rho_c(t) = m_c / |\omega_c(0)| > 0$ is necessarily positive for all times.  Since the initial entropy is non-negative and the entropy in each element is non-decreasing due to $\Pi_{pc} \geq 0$ in \eqref{eqn.euc_S} and due to the absence of source terms in \eqref{eqn.sdthc_S} we have $S_c \geq 0$ for all times. As an immediate consequence, the pressure and the temperature in each cell are always positive, $p_c \geq 0$ and $\theta_c \geq 0$ according to their definition \eqref{eqn.en_pot} and the assumed convexity of the equation of state. 
\end{proof}	

\section{Numerical results} \label{sec.results}
We present a suite of classical numerical test cases for Lagrangian hydrodynamics that aim at showing the correct convergence and shock capturing properties of the novel hyperbolic and thermodynamically compatible Lagrangian scheme, which will be referred to as HTC-Lag method in the following. Differently from most of the existing cell-centered and staggered Lagrangian numerical methods that discretize the Lagrangian hydrodynamics equations in the form \eqref{eqn.PDE_energy}, the energy equation is obtained as a \textit{consequence} of the thermodynamically compatible discretization of the governing PDE system \eqref{eqn.PDE_entropy}, since we are directly solving the \textit{entropy inequality}. Consequently, the major objective of the following test problems is to demonstrate the capability of the HTC-Lag schemes to compute correct solutions for problems with shock waves. 
Depending on the test case, we either use the \textit{a priori} or the \textit{a posteriori} troubled cell detector to compute the blending factor introduced in Section \ref{ssec.Hybridization}. If not stated otherwise, the CFL number is set to $\textnormal{CFL}=0.4$ and the polytropic index of the gas is assumed to be $\gamma=7/5$. For comparison, some of the results are also obtained with the classical cell-centered Lagrangian Finite Volume method,
named EUCCLHYD \cite{EUCCLHYD}, that constitutes the main building block for the dissipative part of our novel HTC-Lag method.  

For all test cases, time integration is performed using the classical explicit fourth order Runge-Kutta method \cite{LeVeque:2002a}, thus ensuring that the GCL and total energy conservation are satisfied up to sufficient accuracy at the fully-discrete level. 
\subsection{Numerical convergence study} \label{ssec.conv}
To verify the accuracy of the new HTC-Lag scheme we rely on the isentropic vortex test case proposed in \cite{HuShuTri}, which is defined in the computational domain $(x,y) \in [0;10]\times [0;10]$ with periodic boundaries everywhere. The initial condition is given in terms of a homogeneous background field $(\rho,\vv,p)^0$ that is supplemented with a set of perturbations defined in \cite{HuShuTri}. The simulations are run until the final time $t_f=1$ on a set of refined meshes with characteristic mesh size $h=\max \limits_c \sqrt{|\Omega_c|}$, and the convergence rates are reported in Table \ref{tab.convRates} using a convective background velocity of the vortex $\vv=(1,1)$. No troubled cell detector is used since the flow field is smooth. The errors are measured in $L_2$ norm for density, horizontal velocity and total energy density, which is obtained via a \textit{post-processing} of the state variables. One can observe that the expected first order of accuracy is achieved.

\begin{table}[!htbp]  
	\caption{Numerical convergence results for the isentropic vortex problem using the HTC-Lag scheme. The errors are measured in the $L_2$ norm and refer to the variables $\rho$ (density), $u$ (horizontal velocity) and total energy $\ee$ at time $t_{f}=1$.}  
	\begin{center} 
		\begin{small}
			\renewcommand{\arraystretch}{1.2}
			\begin{tabular}{c|cccccc}
				$h$ & $||\rho||_2$ & $O(\rho)$ & $||u||_2$ & $O(u)$ & $||\ee||_2$ & $O(\ee)$ \\ 
				\hline
				3.254E-01 & 2.707E-01 & -    & 1.472E-01 & -    & 3.207E-01 & -    \\
				2.490E-01 & 2.171E-01 & 0.82 & 1.125E-01 & 1.00 & 2.544E-01 & 0.86 \\
				1.654E-01 & 1.485E-01 & 0.93 & 7.608E-02 & 0.96 & 1.719E-01 & 0.96 \\
				1.283E-01 & 1.115E-01 & 1.13 & 5.958E-02 & 0.96 & 1.368E-01 & 0.90 \\
			\end{tabular}
		\end{small}
	\end{center}
	\label{tab.convRates}
\end{table} 

\subsection{Riemann problems} \label{ssec.RP}
\begin{figure}[!htbp]
	\begin{center}
		\begin{tabular}{cc}        
			\includegraphics[trim=10 10 10 10,clip,width=0.47\textwidth]{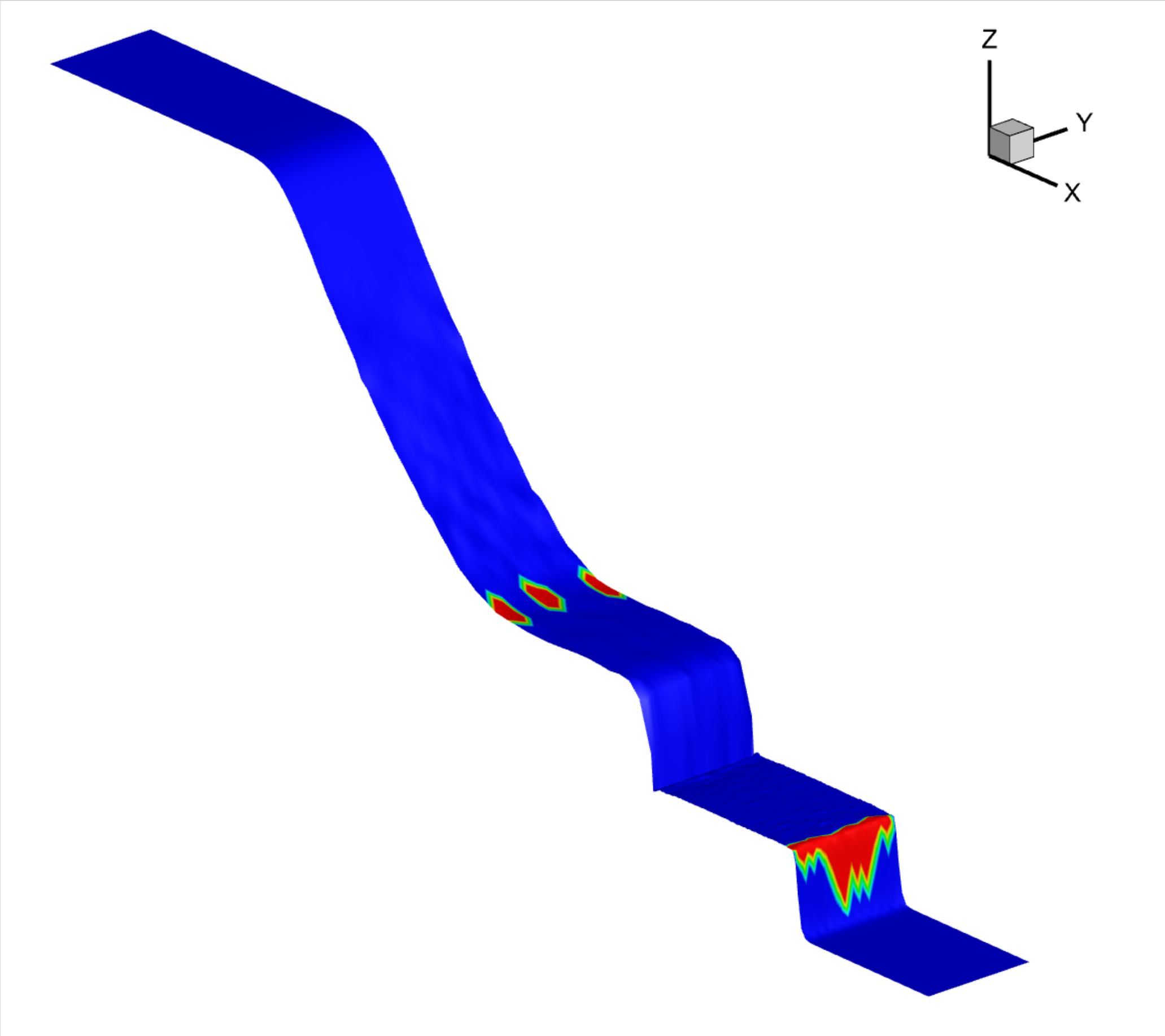} &
			\includegraphics[width=0.47\textwidth]{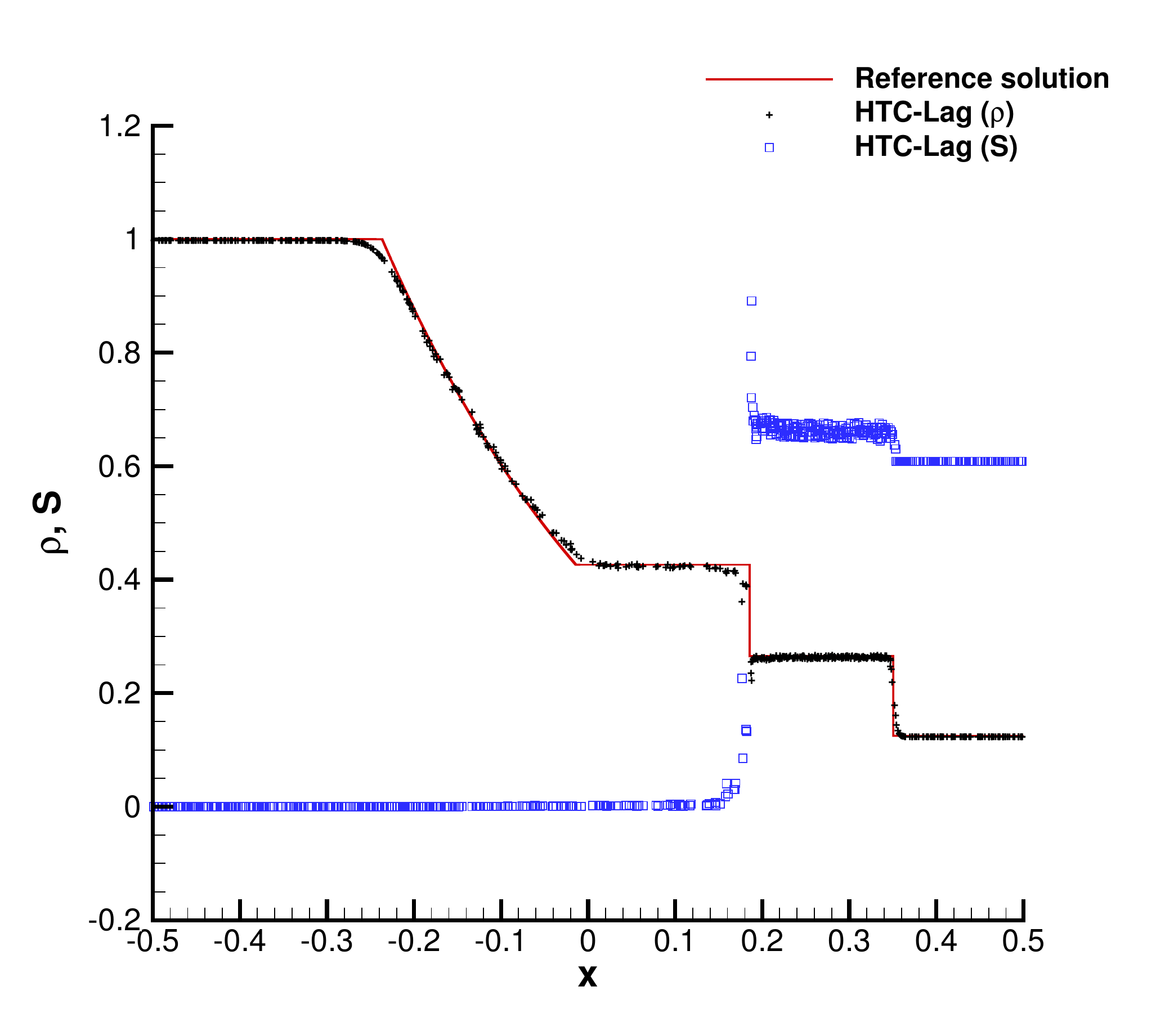} \\
			\includegraphics[width=0.47\textwidth]{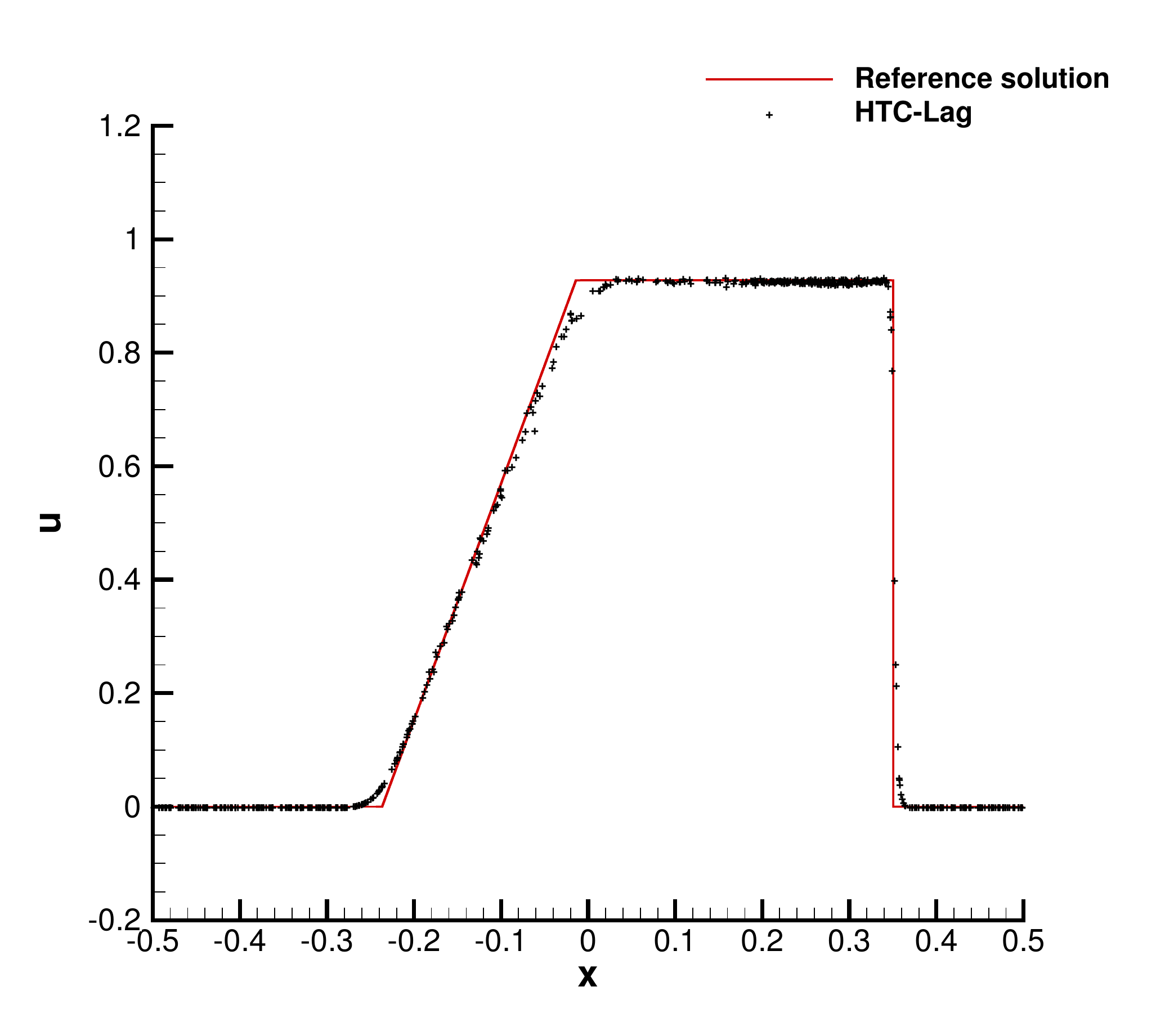} &
			\includegraphics[width=0.47\textwidth]{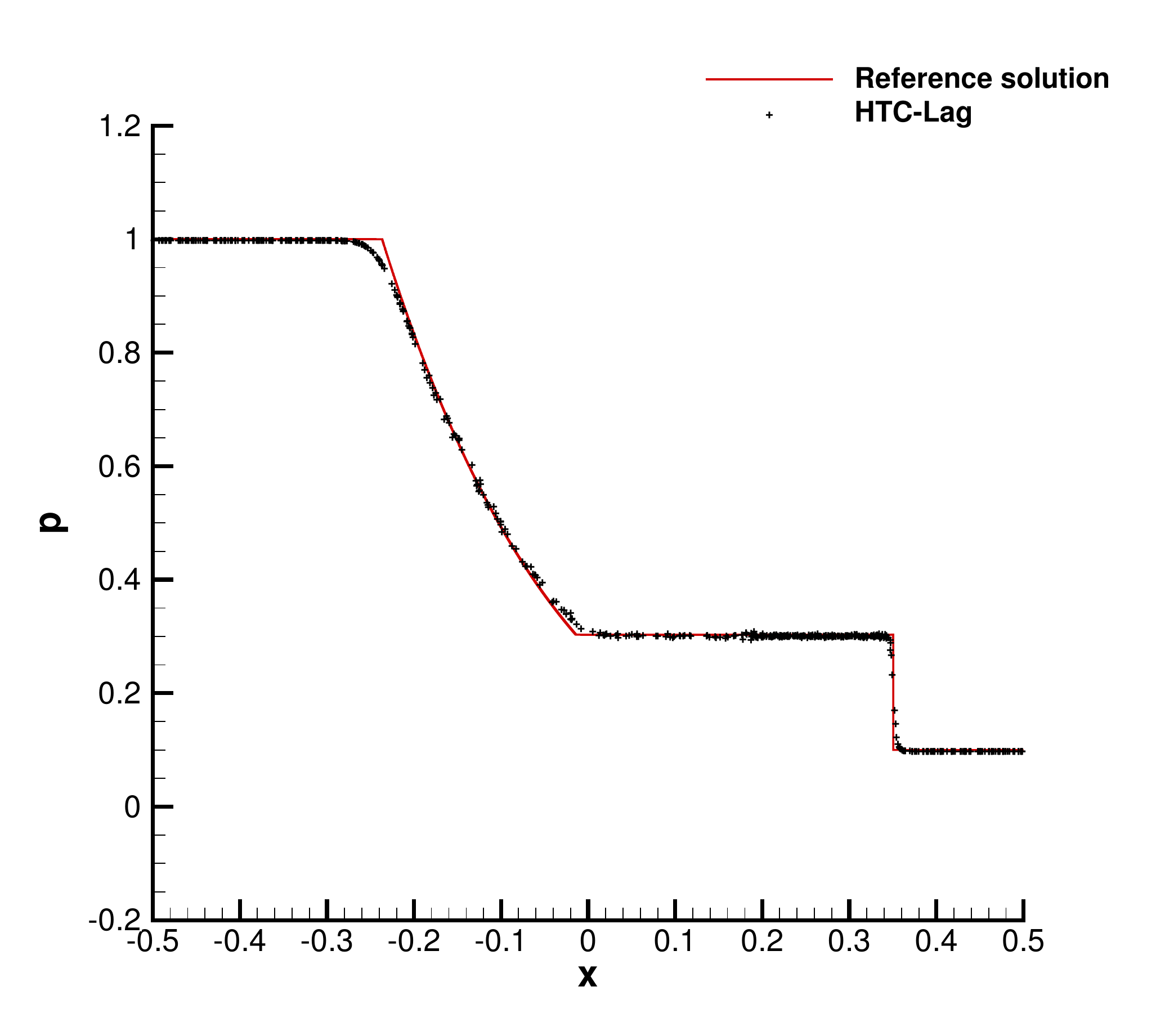} \\
		\end{tabular} 
		\caption{Riemann problem RP1 at time $t=0.2$. Top: detector (left) and scatter plot of density and entropy distribution (right). Bottom: scatter plot of horizontal velocity (left) and pressure (right)}
		\label{fig.Sod}
	\end{center}
\end{figure}

\begin{figure}[!htbp]
	\begin{center}
		\begin{tabular}{cc}        
			\includegraphics[trim=10 10 10 10,clip,,width=0.47\textwidth]{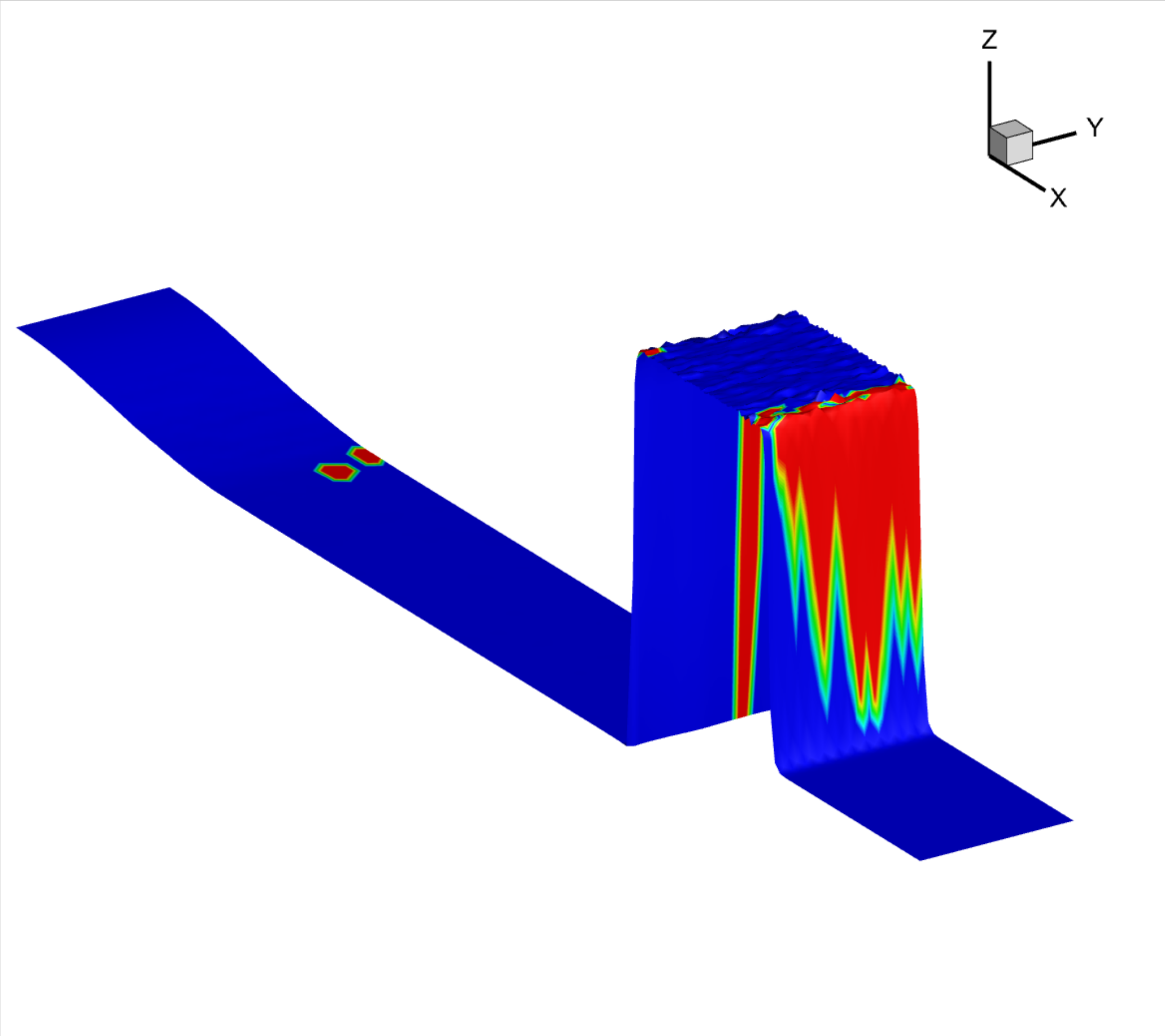} &
			\includegraphics[width=0.47\textwidth]{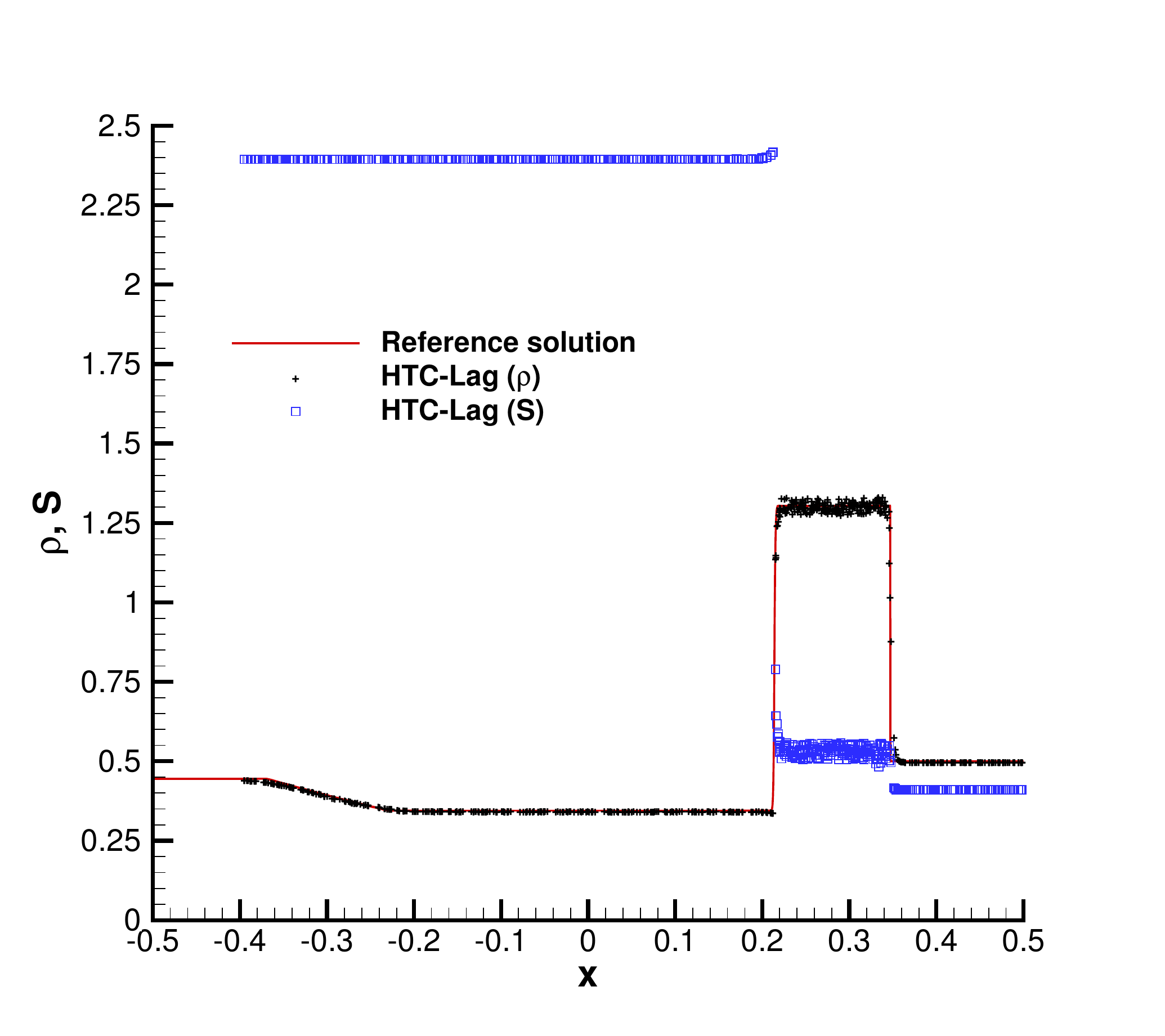} \\
			\includegraphics[width=0.47\textwidth]{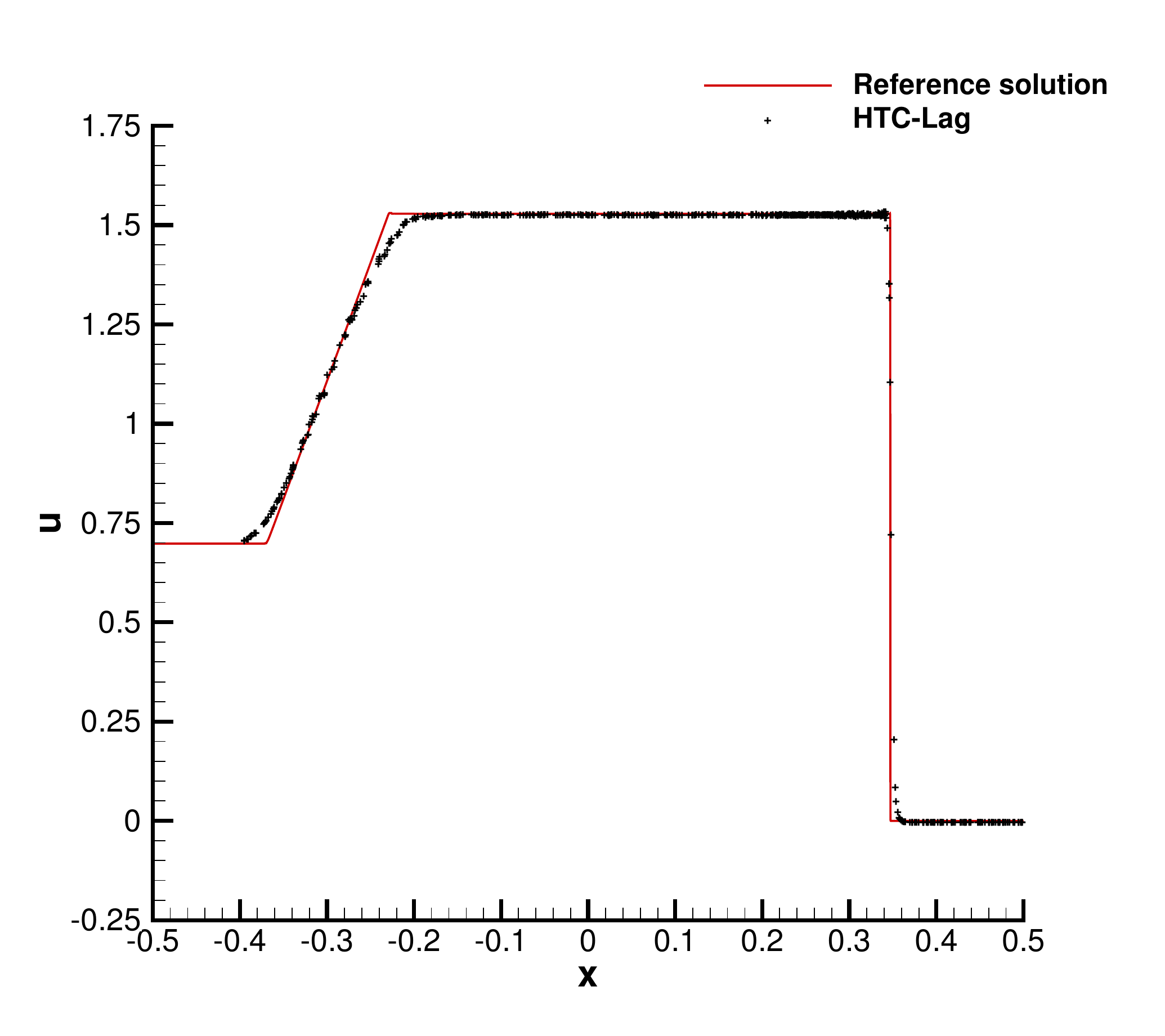} &
			\includegraphics[width=0.47\textwidth]{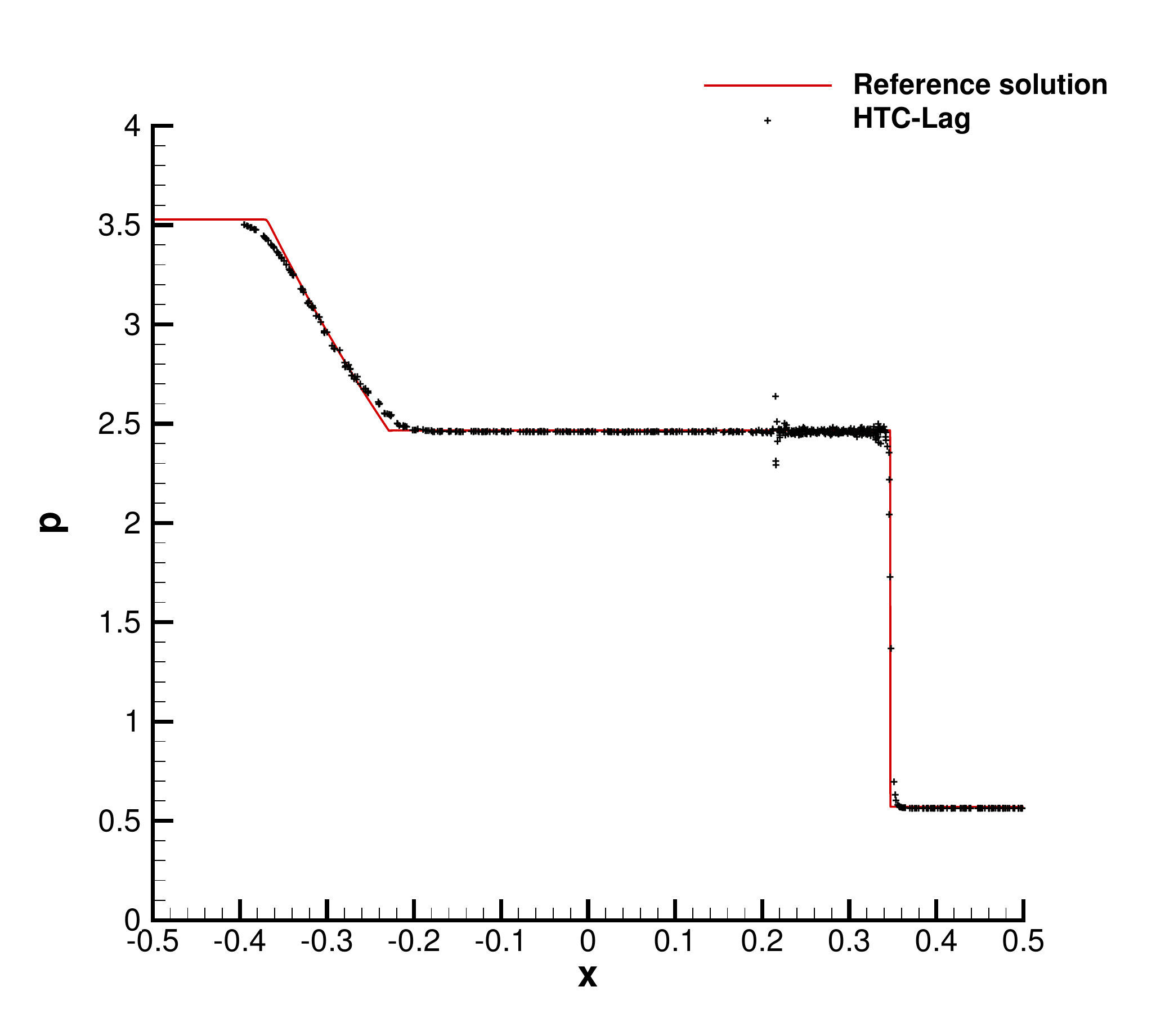} \\
		\end{tabular} 
		\caption{Riemann problem RP2 at time $t=0.14$. Top: detector (left) and scatter plot of density and entropy distribution (right). Bottom: scatter plot of horizontal velocity (left) and pressure (right)}
		\label{fig.Lax}
	\end{center}
\end{figure}

\begin{figure}[!htbp]
	\begin{center}
		\begin{tabular}{cc}        
			\includegraphics[width=0.47\textwidth]{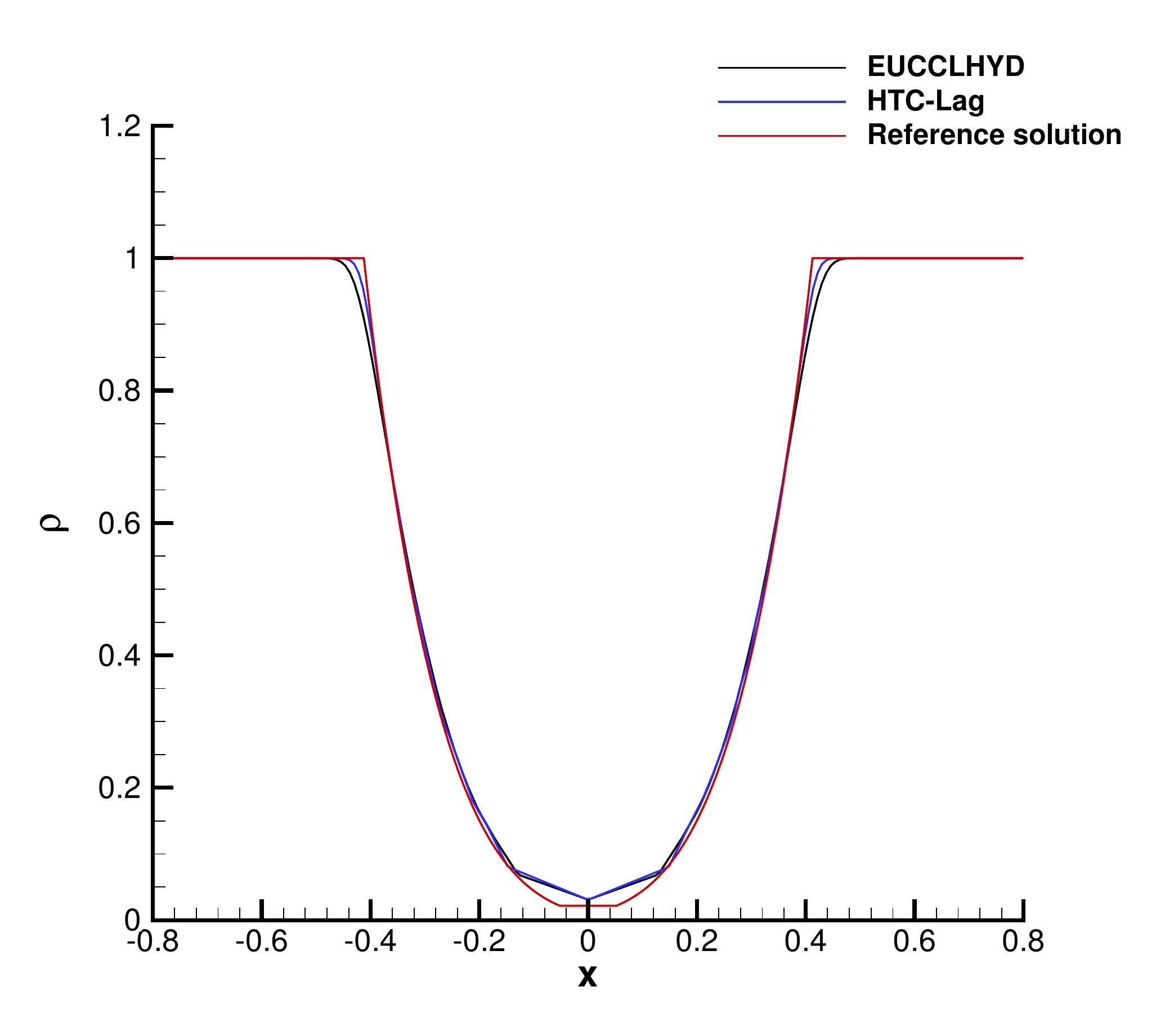} &
			\includegraphics[width=0.47\textwidth]{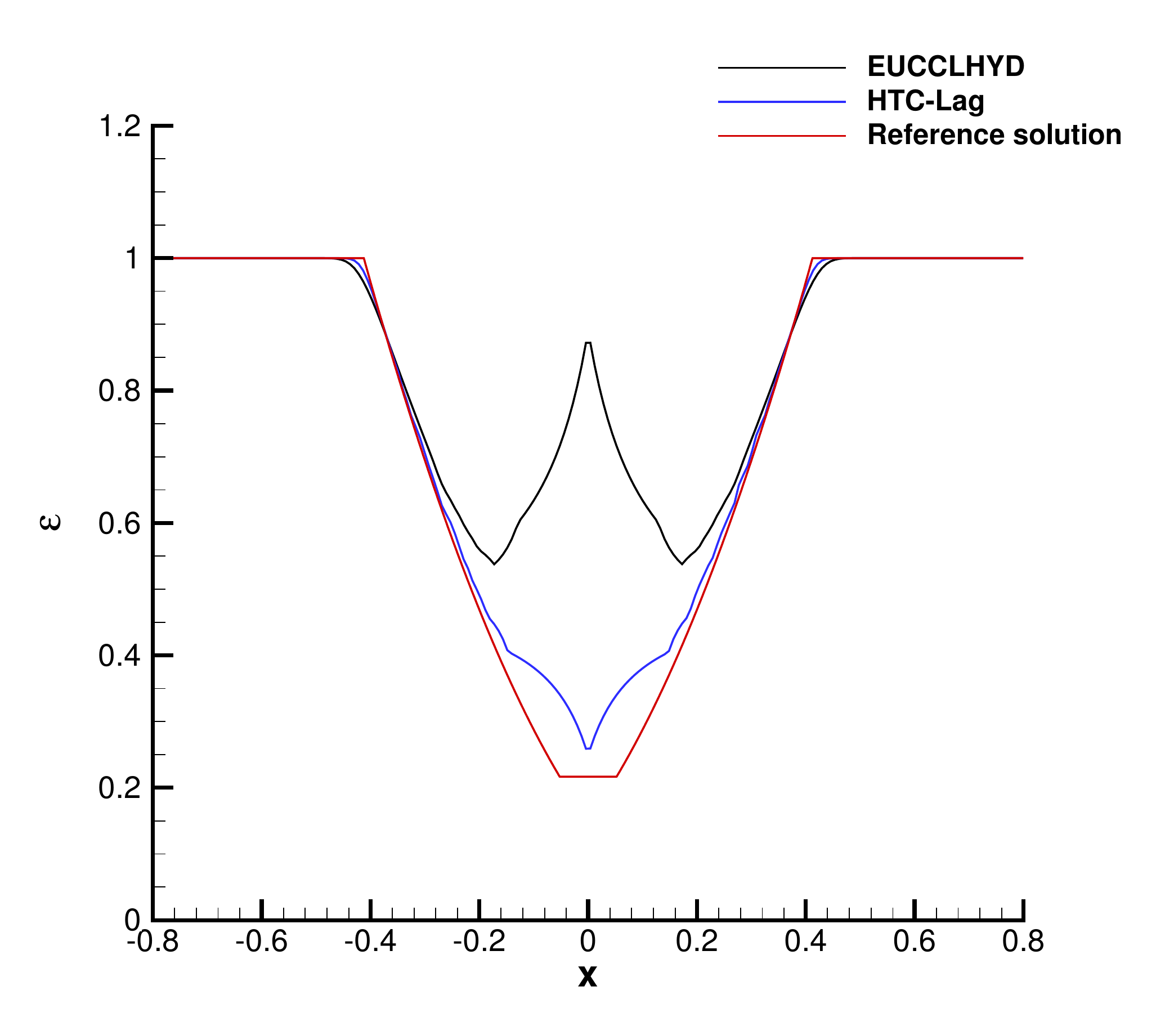} \\
		\end{tabular} 
		\caption{Riemann problem RP3 at time $t=0.15$. Comparison between EUCCLHYD and HTC-Lag results for density (left) and specific internal energy (right).}
		\label{fig.RP123}
	\end{center}
\end{figure}

The novel HTC-Lag scheme is now tested against a set of Riemann problems (RP) with initial data given in Table \ref{tab.RPIC}. The left and right state are given by $(\rho,u,p)_{L,R}$, the initial discontinuity is located at $x_d=0$ and the computational domain is given by $(x,y)\in[x_L;x_R] \times [-0.05;0.05]$ and we use a two-dimensional computational grid with characteristic mesh size of $h=1/200$. Periodic boundary conditions are set in $y-$direction, while the $x-$boundaries are treated as no-slip walls. Despite the one-dimensional setup of these test problems, the simulations are truly multidimensional on unstructured grids, since no element edges are in principle aligned with the main flow. As such, this suite of Riemann problems also allows to check the symmetry preservation of the HTC-Lag method. The numerical solutions are compared with the exact solution of the Riemann problems provided in \cite{ToroBook}. RP1 and RP2 are the Sod and the Lax shock tube problem, respectively, while RP3 corresponds to the double rarefaction test case forwarded in \cite{ToroBook}. We use the \textit{a posteriori} detector of troubled cells for RP1 and RP2, that supplements the scheme with numerical dissipation only across the shock waves, as clearly visible in the first panels of Figures \ref{fig.Sod}-\ref{fig.Lax}. We also note that the results for horizontal velocity and pressure distribution are in excellent agreement with the exact solution, hence proving the correct entropy production by the HTC-Lag scheme across shock waves and to properly resolve rarefaction fans without spurious entropy production. Figure \ref{fig.RP123} shows a comparison between the EUCCLHYD and the HTC-Lag scheme in terms of specific internal energy for RP3, which consists in an isentropic expansion of the gas. It is evident that the new methods are able to properly capture the correct behavior of the gas close to vacuum in the vicinity of the center of the domain. On the other hand, the EUCCLHYD scheme suffers from a spurious entropy production in the rarefaction, which implies a spurious heating.

\begin{table}[!htbp]  
	\caption{Initial states left (L) and right (R) for density, horizontal velocity and pressure for a set of Riemann problems solved on the domain $(x,y)\in[x_L;x_R]$. The final time of the simulation $t_f$ is also indicated.}  
	\begin{center} 
		\begin{small}
			\renewcommand{\arraystretch}{1.2}
			\begin{tabular}{c|cccccc|c}
				\hline
				RP & $\rho_L$ & $u_L$ & $p_L$ & $\rho_R$ & $u_R$ & $p_R$ & $t_f$ \\ 
				\hline
				RP1 & 1.0\phantom{00} & \phantom{-}0.0\phantom{00} & 1.0\phantom{00} & 0.125 & 0.0\phantom{00} & 0.1\phantom{00} & 0.20\\
				RP2 & 0.445 & \phantom{-}0.698 & 3.528 & 0.5\phantom{00} & 0.0\phantom{00} & 0.571 & 0.14\\
				RP3 & 1.0\phantom{00} & -2.0\phantom{00} & 0.4\phantom{00} & 1.0\phantom{00} & 2.0\phantom{00} & 0.4\phantom{00} & 0.15\\
			\end{tabular}
		\end{small}
	\end{center}
	\label{tab.RPIC}
\end{table} 

\subsection{Sedov problem} \label{ssec.Sedov}
\begin{figure}[!htbp]
	\begin{center}
		\begin{tabular}{cc} 
			\multicolumn{2}{c}{		\includegraphics[trim=10 10 10 10,clip,width=0.99\textwidth]{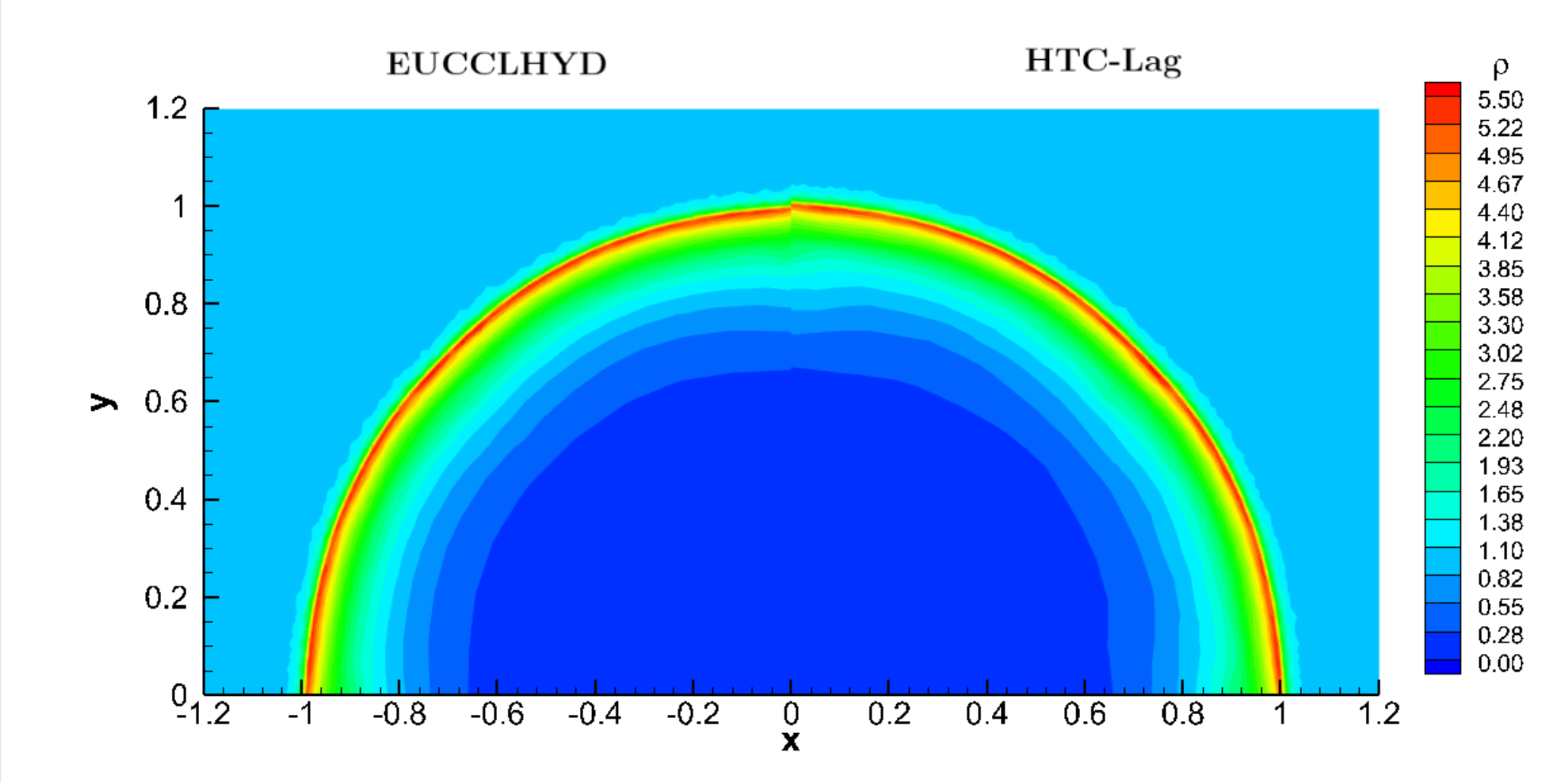}} \\           
			\includegraphics[width=0.47\textwidth]{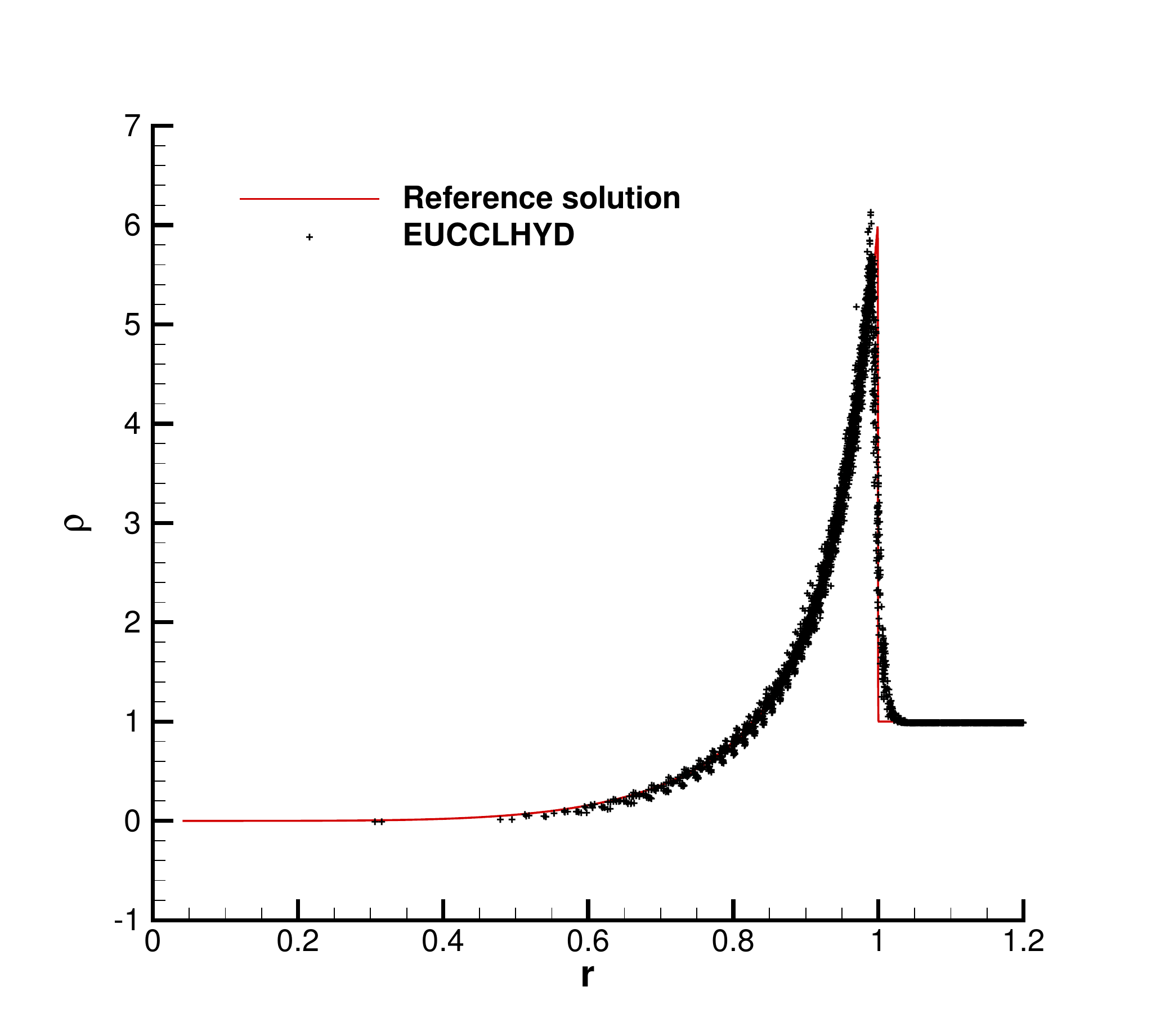} &
			\includegraphics[width=0.47\textwidth]{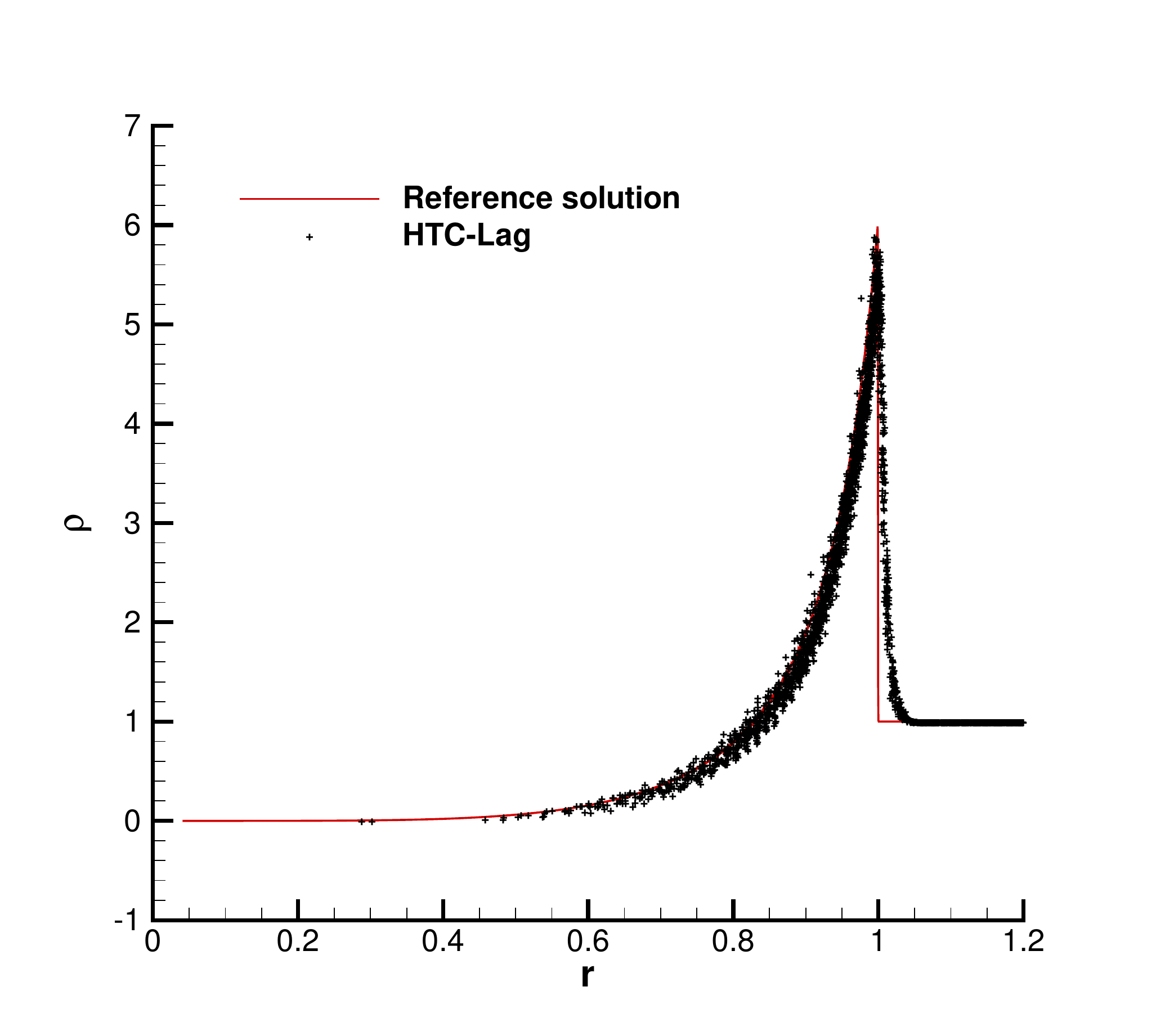} \\
		\end{tabular} 
		\caption{Sedov problem at time $t=1$. Comparison against the EUCCLHYD scheme (left) and the HTC-Lag scheme (right). Density distribution (top) and radial scatter plot of density with the reference solution (bottom).}
		\label{fig.Sedov}
	\end{center}
\end{figure}

\begin{figure}[!htbp]
	\begin{center}
		\begin{tabular}{cc}        
			\includegraphics[width=0.47\textwidth]{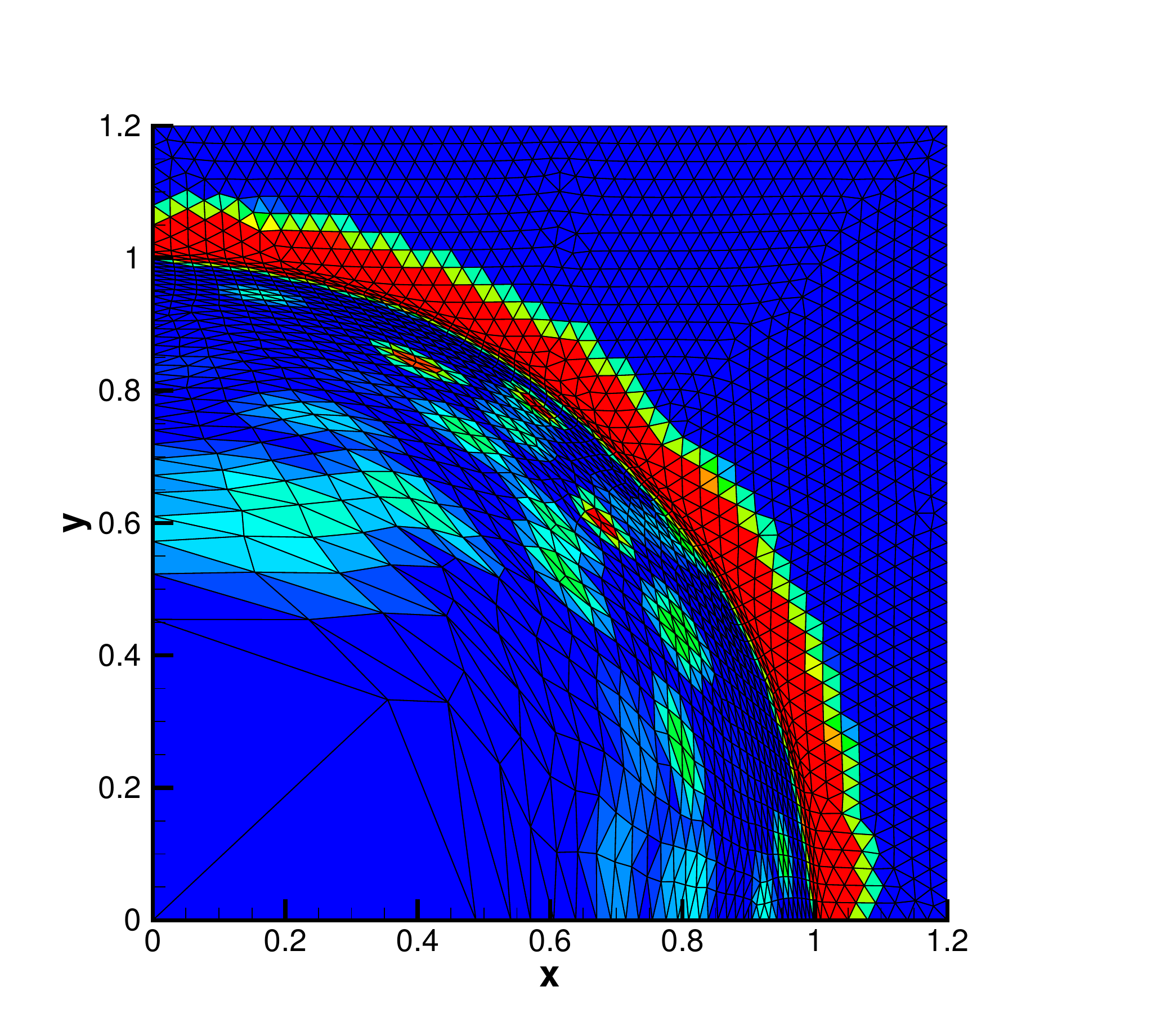} &
			\includegraphics[width=0.47\textwidth]{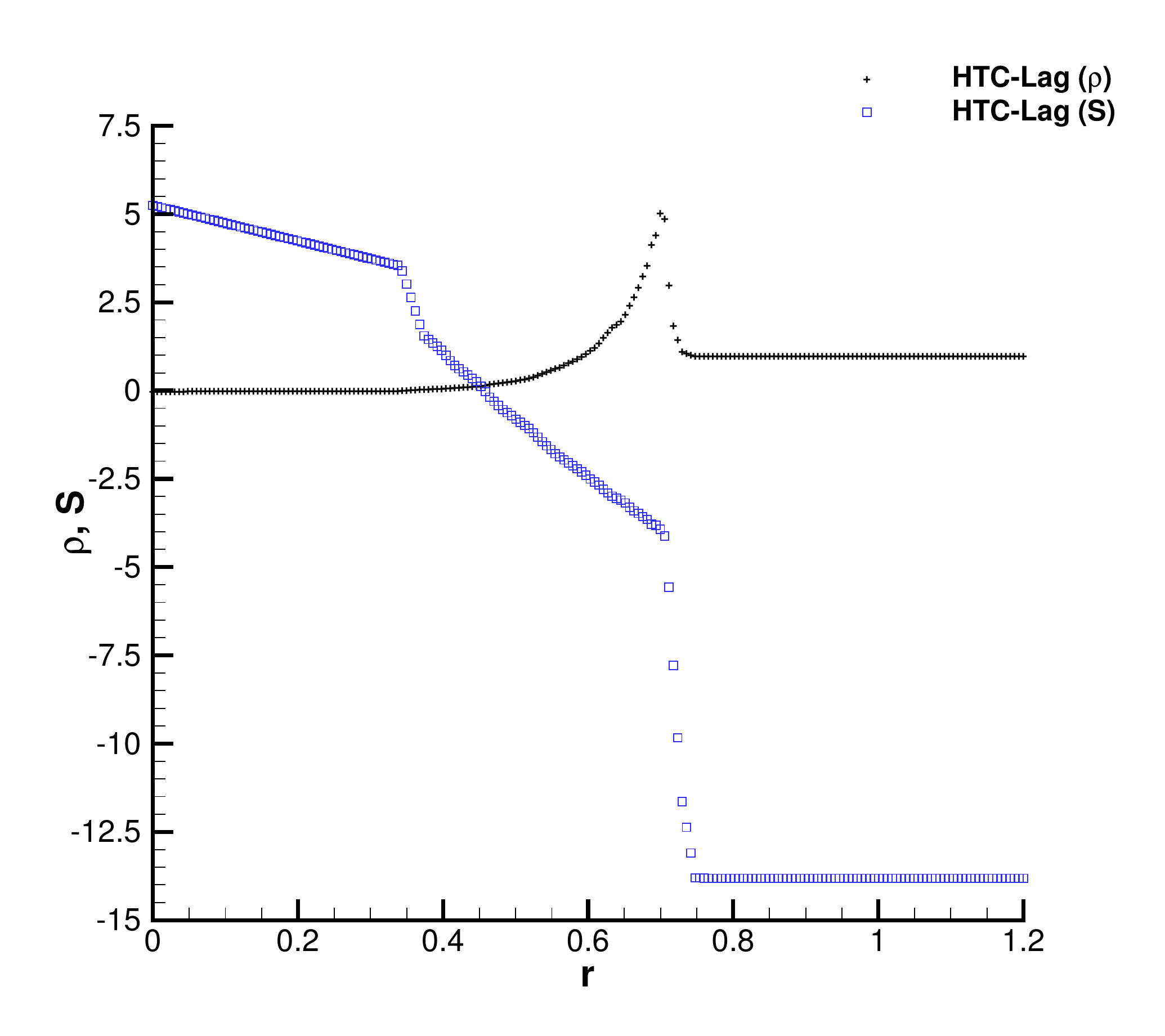} \\
		\end{tabular} 
		\caption{Sedov problem at time $t=1$. Detector (left) and radial density and entropy distribution (right).}
		\label{fig.Sedov_rho_S}
	\end{center}
\end{figure}

This problem involves a strong blast wave departing from the origin at $(0,0)_{or}$ of the computational domain $(x,y) \in [0;1.2]\times [0;1.2]$. We set wall boundary conditions on the left and bottom side, while transmissive pressure boundaries are imposed elsewhere. The domain is paved with an unstructured triangular mesh with characteristic size of $h=0.03$ yielding a total number of $N_c=3602$ cells. The fluid is initially assigned with a constant state given by $(\rho,\vv,p)^0=(1,\mathbf{0},10^{-6})$. The energy of the explosion is totally concentrated in the cells containing the origin by prescribing an initial pressure $p_{or}=(\gamma-1)\rho^0 \varepsilon^0/|\omega|_{or}$, where $|\omega|_{or}$ denotes the volume of the cells containing the origin. According to \cite{SedovExact}, the amount of released energy is set to $\varepsilon^0=0.244816$, thus the solution consists of a diverging infinite strength shock wave that is placed at radius $r=1$ at the final time of the simulation $t_f=1$. The analytical solution can be found in \cite{SedovExact}.   
This test is run using both the new HTC-Lag and the EUCCLHYD scheme, and the final density distribution is qualitatively compared in Figure \ref{fig.Sedov}. Being less dissipative, the HTC-Lag scheme better captures the shock peak and location, which is also visible from the density scatter plots. We use the \textit{a priori} detector of troubled cells that is only active across the shock wave, as expected, thus entropy increases perfectly in correspondence of the discontinuity, as shown in Figure \ref{fig.Sedov_rho_S}.

\subsection{Cylindrical expansion into vacuum} \label{ssec.VacuumExp}
\begin{figure}[!htbp]
	\begin{center}
		\begin{tabular}{cc}        
			\includegraphics[trim=10 10 10 10,clip,width=0.47\textwidth]{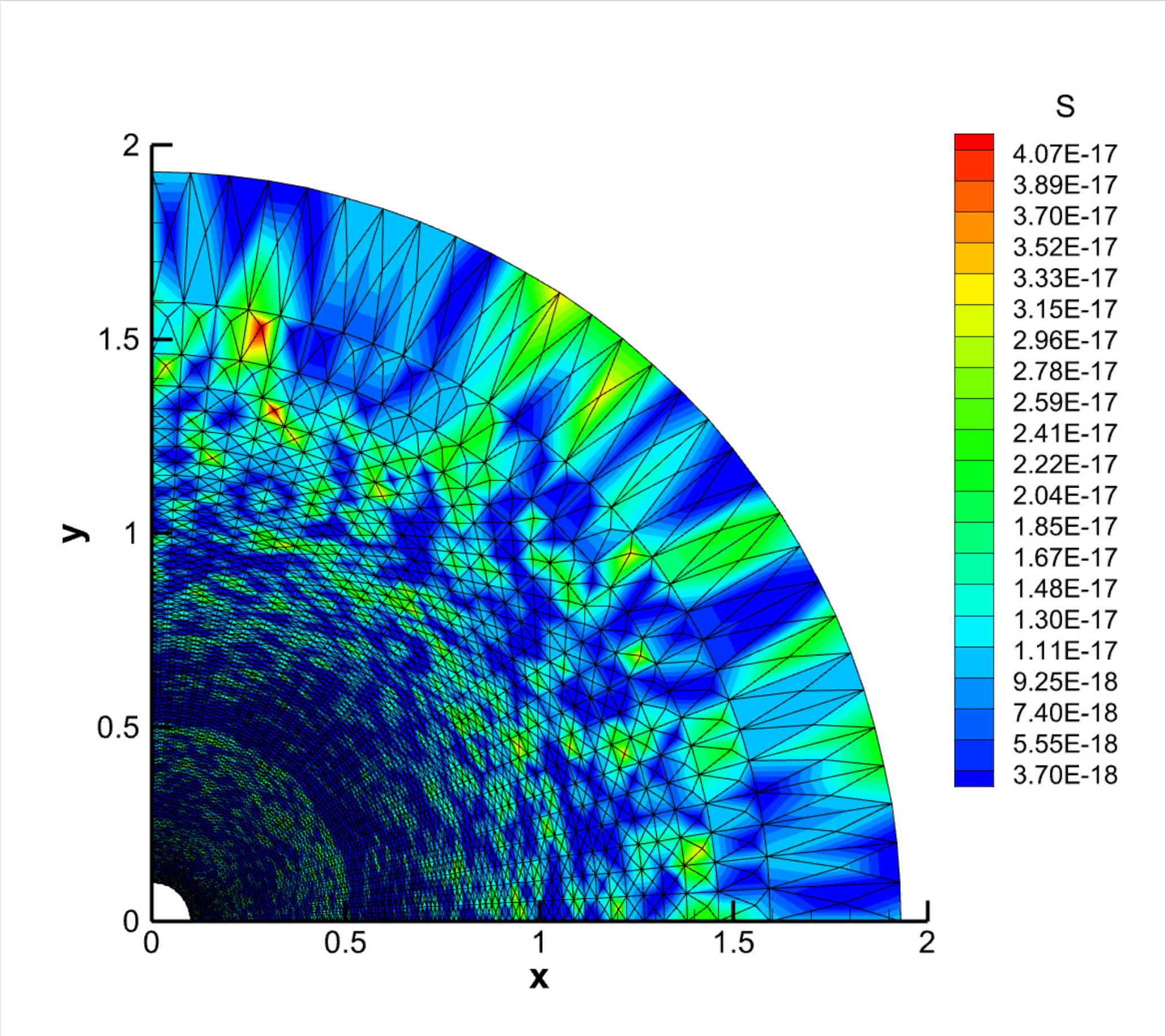} &
			\includegraphics[trim=10 10 10 10,clip,width=0.47\textwidth]{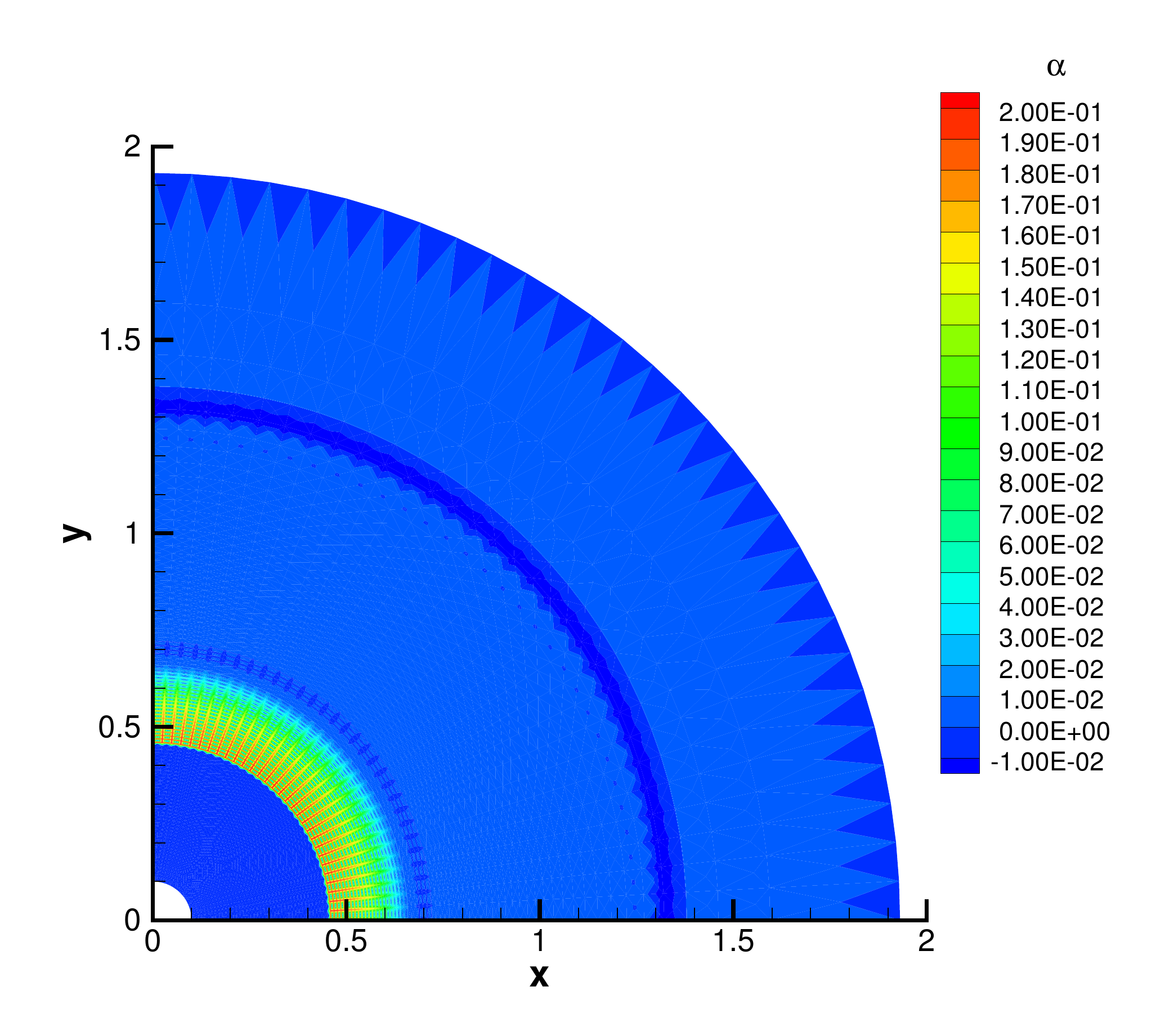}\\
        	\includegraphics[width=0.47\textwidth]{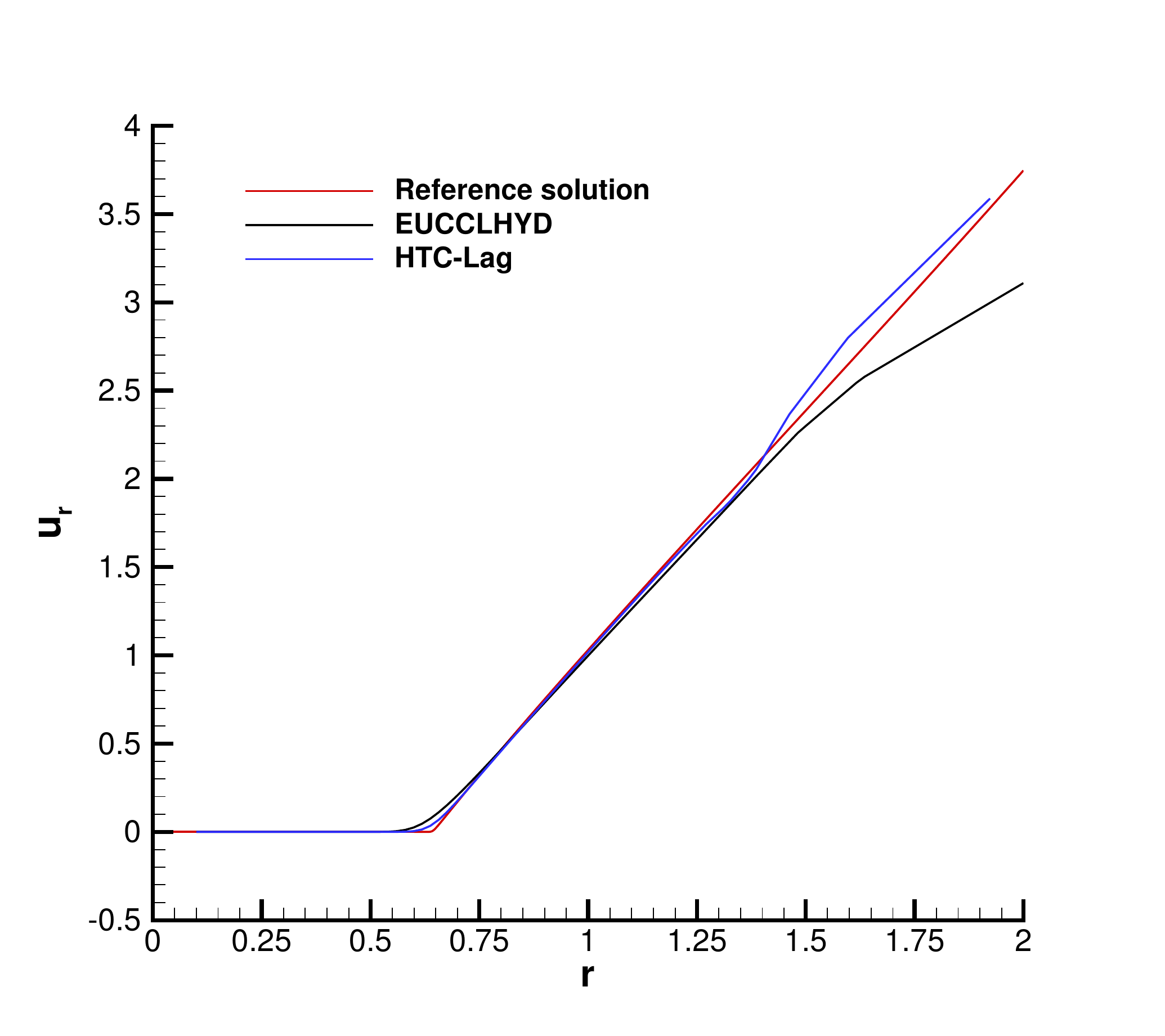} &			
			\includegraphics[width=0.47\textwidth]{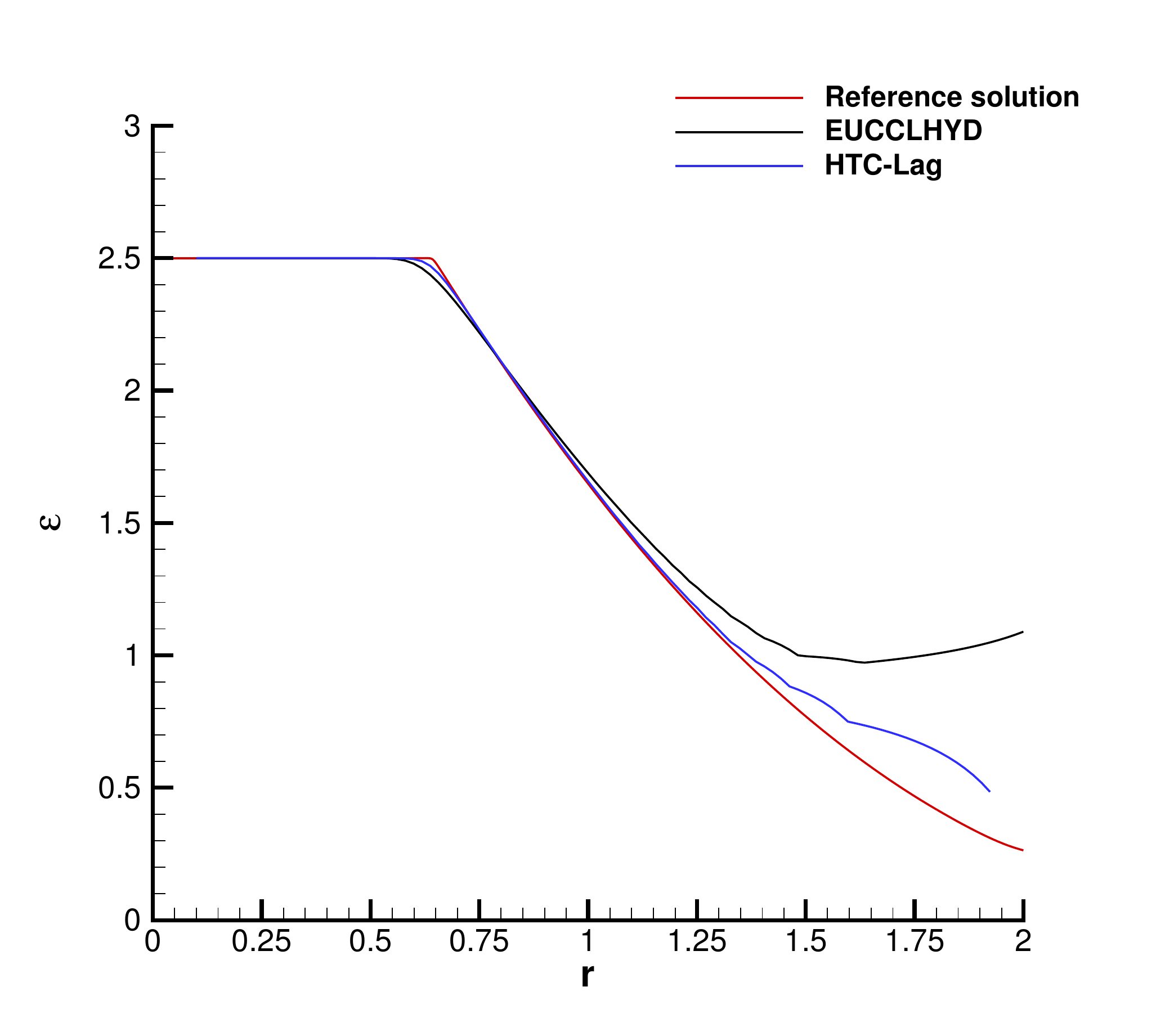} \\
		\end{tabular} 
		\caption{Cylindrical expansion into vacuum $t=0.3$. Top: entropy distribution with mesh configuration (left) and distribution of the correction scalar $\alpha$ (right). Bottom: comparison of the HTC-Lag scheme with EUCCLHYD scheme and reference solution for radial velocity (left) and specific internal energy (right).}
		\label{fig.VacuumExp}
	\end{center}
\end{figure}

As last test case we solve a cylindrical expansion of gas into vacuum, as proposed in \cite{Maire2020}. The computational domain is given by the portion of a shell defined with cylindrical coordinates $(r,\phi) \in [0.1;0.1]\times [0;\pi/2]$, and it is paved with a triangular grid of sizes $h_r=1/100$ and $h_{\phi}=1/30$ in radial and angular direction, respectively. Wall boundaries are prescribed at $\phi=0$ and $\phi=\pi/2$, whereas we apply a free pressure ($p_b=0$) boundary condition on the internal and the external frontier of the shell. No troubled cell detection is carried out for this test case, thus the HTC-Lag scheme is run with zero entropy production until the final time $t_f=0.3$. Figure \ref{fig.VacuumExp} shows the final mesh configuration with the entropy distribution, which remains of machine accuracy throughout the entire computation. We also plot the scalar  correction factor $\alpha_p$, which is mostly active across the density and pressure jumps, detected with the \textit{a posteriori} indicator. Finally, the results are compared against those obtained with the EUCCLHYD scheme in terms of radial velocity and specific internal energy distribution. Furthermore, a reference solution is obtained by solving this test case on a very fine one-dimensional mesh with 20'000 cells by a second order MUSCL-Hancock-type TVD scheme. We observe that the classical cell-centered formulation is plagued by a spurious heating in the vicinity of the interface with the vacuum. The new thermodynamically compatible schemes can cure quite well this behavior, obtaining results which are in very good agreement with the reference solution, even in the vicinity of the vacuum interface.

\section{Conclusions} \label{sec.concl}
In this work a new cell-centered Lagrangian finite volume scheme for gas dynamics has been designed that is by construction compatible with the second law of thermodynamics and which satisfies the total energy conservation as a consequence at the semi-discrete level. Unlike existing schemes, the novel method \textit{directly} evolves the \textit{entropy} as a primal evolution quantity while total energy conservation is merely retrieved as a \textit{consequence} of the new thermodynamically compatible scheme. The new method is proven to be nonlinearly stable in the energy norm. Moreover, we demonstrate the positivity preserving property for both density and pressure. Entropy conservative as well as entropy stable schemes are derived, and they are blended by a space-time adaptive technique which preserves all the thermodynamic properties of the original methods. Numerical convergence  rates and results for classical benchmarks in Lagrangian hydrodynamics have been run, showing the advantages of the new thermodynamically compatible schemes over standard cell-centered Lagrangian methods, especially in regions of isentropic flow. 

In the future we plan to design a thermodynamically compatible Lagrangian discretization for the GPR model of continuum mechanics \cite{PeshRom2014}, starting from the Lagrangian finite volume scheme proposed in \cite{LGPR}. We also aim at investigating the development of a more sophisticated blending factor in the adaptive scheme, in order to provide a smoother transition between entropy conservative/stable schemes.

\section*{Acknowledgments}
WB received financial support by Fondazione Cariplo and Fondazione CDP (Italy) under the grant No. 2022-1895 and by the Italian Ministry of University and Research (MUR) with the PRIN Project 2022 No. 2022N9BM3N.
MD was funded by the Italian Ministry of University and Research (MUR) in the frame of the PRIN 2022 project \textit{High order structure-preserving semi-implicit schemes for hyperbolic equations} and by the European Union NextGenerationEU (PNRR, Spoke 7 CN HPC). Views and opinions expressed are however those of the author(s) only and do not necessarily reflect those of the European Union or the European Research Council. Neither the European Union nor the granting authority can be held responsible for them.

\bibliographystyle{plain}
\bibliography{biblio}
\end{document}